\documentclass[11pt,twoside,leqno]{article}
\usepackage{epsfig}
\newenvironment{proof}{\noindent {\bf Proof.}\hspace{2mm}}{\qed}

\DeclareRobustCommand{\qed}{%
  \ifmmode \mathqed
  \else
    \leavevmode\unskip\penalty9999 \hbox{}\nobreak\hfill
    \quad\hbox{$\square$}\medbreak%
  \fi
}
\newcommand{\mathqed}{\quad\hbox{$\square$}}
\newtheorem{lemma}{Lemma}[section]
\newtheorem{theorem}{Theorem}
\newtheorem{proposition}[lemma]{Proposition}

\newtheorem{remark}[lemma]{Remark}
\newtheorem{corollary}[lemma]{Corollary}
\newtheorem{construction}{Construction}
\newtheorem{definition}[lemma]{Definition}

\newtheorem{conjecture}{Conjecture}

\newcommand\acl{\hbox{\rm acl}}

\newcommand\implies{\Rightarrow}
\newcommand\includedin{\subseteq}

\newcommand\union{\cup}

\newfont{\Bbbb}{msbm10}

\newcommand\Zz{\mbox {\Bbbb Z}}

\newcommand\CC{{\cal C}}

\input amssym.def
\input amssym

\newsymbol\square 1003
\newsymbol\rtimes 226F
\newsymbol\restrictedto 1316

\makeatletter
\newcommand\myeqnlabel[1]{%
  \hb@xt@ .01\p@ {}\rlap {\normalfont \normalcolor 
  \hskip -\displaywidth {#1}}}
\newenvironment{texteqn}[2][0.8]
    {\begin{equation}\renewcommand\@eqnnum{\myeqnlabel{#2}}
     \begin{minipage}{#1\linewidth}}
    {\end{minipage}\end{equation}}

\newtheorem{rawnamedtheorem}{}

\newenvironment{namedtheorem}[1]{\renewcommand{\therawnamedtheorem}{#1}
   \begin{rawnamedtheorem}}
  {\end{rawnamedtheorem}}
\makeatother 
\begin{document}
\title{Universal graphs with a forbidden subtree}
\author{Gregory Cherlin\thanks{First author supported by NSF Grant DMS 0100794.}\\
Department of Mathematics, Rutgers University\\
110 Frelinghuysen Rd., Hill Center\\
 Piscataway, New Jersey 08854, U.S.A\\
\and
Saharon Shelah\thanks{Second author's Research supported in part by U.S.-Israel
Binational Science Foundation grant 0377215; collaboration 
   supported in part by NSF Grant DMS-0100794. Paper 850.}\\
Dept. Mathematics, Hebrew University, Jerusalem, Israel\\
and Rutgers University, Piscataway, NJ 08854
}

\maketitle

\section{Introduction}

The systematic investigation of countable universal
graphs with ``forbidden'' subgraphs was initiated in \cite{kp1},
followed by \cite{kmp}.  If $C$ is a finite connected graph, then a
graph $G$ is {\it $C$-free} if it contains no subgraph isomorphic to
$C$.  A countable $C$-free graph $G$ is {\it weakly universal} if
every countable $C$-free graph is isomorphic to a subgraph of $G$, and
{\it strongly universal} if every such graph is isomorphic to an {\it
induced} subgraph of $G$.  Such universal graphs, in either sense, are
rare.  
Graph theorists tend to use the term ``universal'' in
the weak sense, while model theorists tend to use it in the strong
sense. We will use the term in the graph-theoretical sense here:
``universal'' means ``weakly universal'', though we sometimes
include the adverb for emphasis.

We deal here with the problem of determining the 
finite connected constraint graphs $C$ for which there is a countable
universal $C$-free graph. We introduce a new inductive method and
use it to settle the case in which $C$ is a tree, confirming a
long-standing conjecture of Tallgren. The existence of such a
countable universal graph says something about the class of all finite
$C$-free graphs, and something about $C$ itself, and the problem
ultimately is to determine what, exactly, it does say. 

This has been partially elucidated in \cite{css}. Associated to any 
constraint graph $C$ there is a natural notion of {\it algebraic closure}: 
loosely speaking, a vertex is in the algebraic closure of a given set
if the number of vertices of the same type must be finite
in any $C$-free graph.
For example, if $C$ is a star consisting of a vertex with a set of
neighbors, so that the $C$-free graphs are those with a fixed bound
on the vertex degrees, the algebraic closure of a set is simply the
union of its connected components. In general, if the algebraic
closure of a set is the union of the algebraic closures of its
elements, we say that the operation is {\it unary}.
This is usually not the case: if for example 
$C$ is a circuit of length $4$, then the associated
algebraic closure operation is generated by a (partial) binary
operation: adjoin the unique common neighbor of any pair of points, if 
it exists; and iterate.
At the opposite extreme, the algebraic closure operation may be
trivial: for example, if $C$ is a complete graph then the algebraic 
closure of a set is the set itself. 

Now it turns out that the following three conditions are intimately
related at both a theoretical and empirical level.
\begin{enumerate}
\item The algebraic closure operation associated to $C$ is {\it
 (uniformly) locally finite} in the sense that the algebraic closure
 of a set of size $n$ in any $C$-free graph is finite (and then
 necessarily bounded in size by a function of $n$).
\item There is a strongly universal $C$-free graph.
\item There is a weakly universal $C$-free graph.
\end{enumerate}

These conditions are successively weaker, and not much different in
practice. We know of no case which separates the second from the third
condition, but there are trivial examples falling under the last two
cases and not the first: for graphs of maximal vertex degree 2 there
is a universal graph made up of infinitely many cycles of all lengths
and infinitely many two-way infinite paths. As the algebraic closure 
of a single vertex is its connected component, local finiteness fails,
but what happens in this case is that the algebraic closure operation
is very tightly structured. 
More generally, the same phenomenon occurs,
and the situation as a whole is very much the same, 
if the constraint graph is a {\it near-path}, that is a tree which
is not a path, but is obtained by attaching one edge with one
additional vertex to a path. 

Our goal is to arrive at a more concrete understanding of the
exceptional constraints allowing a weakly universal $C$-free graph. It
makes good sense to state the problem in full generality as follows:
{\it is there an effective procedure to decide whether a given finite
  constraint $C$ allows a corresponding universal (countable) graph?}

The problem is very much open and has been attacked from two
directions. 
Some encoding results are known which aim in the direction of 
a proof of undecidability.
But the bulk of the research, like our present line,
aims in the opposite direction, developing general tools to settle
problems of this type: indeed, it is plausible at this stage
that it may be possible to work out 
the list of exceptional constraints $C$ allowing a universal 
graph in a completely explicit way. This 
would be the strongest form of a
positive solution to the decision problem, though ``softer''
approaches are also available.
While the condition of local finiteness is essentially a 
{\em halting problem} for a specific computation, the computations 
in question tend in the vast majority of cases to diverge.

The present paper has two goals: to  present a simple
inductive style of argument which suggests that 
if the list of exceptional ``favorable'' constraints 
is in fact as simple as we are suggesting it should be, then 
we should be able to prove that fact, and 
to buttress this claim by establishing Tallgren's Tree Conjecture.
While we do not have a conjectured list of favorable constraints $C$
in general, in the case of trees Tallgren conjectured the simple answer
for this case explicitly many years ago. One might expect that
knowing the answer would be the major ingredient in finding a proof, 
and as far as that goes it probably is, but nonetheless the question
has remained open, and it seems to need our inductive method
in order to be reduced to a finite number of individually treatable
minimal cases.

The result is as follows.

\begin{theorem}[Tree Conjecture]\label{Intro:TC}
If $T$ is a finite tree, then the following are equivalent.
\begin{enumerate}
\item There is a weakly universal $T$-free graph.
\item There is a strongly universal $T$-free graph.
\item $T$ is a path or a near-path.
\end{enumerate}
\end{theorem}

\noindent 
We remark that the algebraic closure operation 
associated to a path is locally finite, and 
as we have indicated near-paths fall into the exceptional class 
behaving much like the star of order four. They may indeed exhaust the
exceptional class for which the associated algebraic closure operation
is not locally finite, but there is an associated universal graph.

While all of this might suggest that only the most obvious
examples of universal graphs can exist, this is not so.
Komj\'ath showed, unexpectedly,
that the 2-bouquet formed by joining two triangles over a common
vertex provides another constraint graph $B$ allowing a (strongly)
universal $B$-free graph \cite{kom}; more generally, 2-bouquets
$B_{m,n}$ formed by joining complete graphs $K_m$ and $K_n$ over a
single common vertex have been
thoroughly analyzed in \cite{ct}: there is a (weakly or strongly)
universal $B_{m,n}$-free graph if and only if the parameters satisfy
the following conditions:
$$\min(m,n)\le 5; \qquad (m,n)\ne (5,5)$$
This is visibly a delicate condition, and requires a close
combinatorial analysis to achieve. Examples of this type continue to
hold open the possibility that the final list may be
more delicate than anything we have seen to date.

There has been prior work on the case of tree constraints.  
First, taking a path or a near-path as forbidden subgraph does allow
a universal graph \cite{kmp,ct}.  In the other direction, the
nonexistence of universal graphs has been treated in the following
cases: (1) arrows, which are trees consisting of a path with two more
edges adjoined to either endpoint, a case treated in \cite{gk}; 
(2) trees with a unique vertex of maximal degree $d\ge 4$, 
which is moreover adjacent to a leaf, treated in \cite{fk2};
and (3) ``bushy'' trees, that is trees
with no vertex of degree 2, treated in \cite{cs1}. Of these, the case
treated in \cite{fk2} now seems the most suggestive. Indeed, we will
show that by combining the case treated in \cite{fk2} with a simple
inductive idea, we deduce that for trees $C$ having a unique vertex of
maximal degree $d\ge 4$, there is no strongly universal $C$-free graph.
The argument is prototypical for the general problem.

What we do here is motivated to a degree by a false conjecture
in \cite{css}, the Monotonicity Conjecture: if a constraint $C$ allows 
a universal graph (in either sense) then for any {\it induced} subgraph
$C_0$ of $C$, the tighter constraint $C_0$ also allows a universal
graph. It turns out that the close analysis of 2-bouquets $B(m,n)$
refutes this, as the graph $B(5,5)$ is an induced subgraph of $B(5,6)$ 
and the latter allows a strongly universal graph while the former does
not allow a weakly universal graph. But there is enough truth to the
conjecture to make it useful: if one passes to an induced subgraph by
the operation of {\it pruning} introduced here (removing certain
$2$-blocks), the monotonicity principle is valid. So after explaining
this, in taking up the Tree Conjecture we will deal with {\it
  critical} trees, which by definition are the trees which are not
paths or near-paths, but become paths or near-paths when
pruned---which, incidentally, is nothing but the removal of leaves in
this case. The reader can see for himself that the structure of these
trees is very simple. What we need to prove is that $(a)$ the algebraic
closure operation is not locally finite in these cases; $(b)$ this
leads to the nonexistence of strongly or weakly universal $T$-free
graphs in all cases. Now $(a)$ is easier than $(b)$ and is a
prerequisite for the latter, and the constructions used to accomplish
$(a)$ serve as templates for the more delicate constructions used to
accomplish $(b)$. In fact there are three layers of constructions.
We found it fairly easy to decorate the constructions used in case
$(a)$ to refute the existence of strongly universal $T$-free graphs
for $T$ a critical tree $T$, and rather troublesome to convert 
the latter into
refutations in the weakly universal case, which is the 
problem which was originally posed. Naturally we suppress all of these
intermediate steps except in some illustrative cases. The result is that
one will see various ``bells and whistles'' in the constructions,
and arguments that to a certain extent the graphs that interest us
are maximal in the sense that an embedding into a larger $T$-free
graph does not create new edges. This is not literally the case: it
would be more accurate to say that certain critical vertices 
acquire no new vertices, and even this overstates the matter.

We make one further remark about these constructions.
With the exception of the first cases treated (called monarchy and
stardom), these constructions do not leap to mind; perhaps with better
insight they should, and in any case as they are all variations on one
theme this theme can be added to the toolbox for future reuse.
But at one point we doubted the truth of the Tree Conjecture,
and computed the algebraic closure operator for the case of a specific
tree on 14 vertices, the ``most likely'' counterexample to the conjecture.
This turned out to be so tightly constrained that it only allowed one
type of construction, which is the one used here throughout.
Similarly, in working with bouquets, 
it is doubtful that one would find the construction used to
refute the existence of a universal $B(5,5)$-graph on an ad hoc basis, 
but again one computes the behavior of the algebraic closure operator
(expecting local finiteness, in fact) and the relevant construction
simply appears. In the case of bouquets, where there are many cases of
local finiteness, one can actually see these computations in
\cite{ct}, where they are necessary for the main results. Here they
can and should be suppressed, as what interests us are the necessary
constructions. But as the paper represents unfinished business, indeed
merely the initial step of what could be a very long process, these
methodological points should be noted.

Furthermore, we have drastically oversimplified our discussion in one
crucial respect. All of these problems make equally good sense---more
sense, in fact---when the constraint $C$ is replaced by a finite set
$\CC$ of finite, connected, constraint graphs. Some things become
clearer in the process: the case of complete constraint graphs
generalizes to the case of sets $\CC$ closed under homomorphism, where
the common feature is that the algebraic closure operation is
{\it trivial}: $\acl(A)=A$ for all $A$.  Furthermore examples of {\it
  mixed type} occur: one may take any constraint $C$ for which the
algebraic closure operator is locally finite, combine it with any
further finite set of constraints closed under homomorphism, without
altering the algebraic closure operator. These phenomena remain
invisible when one considers only single constraints, and a number of
the more general examples are exceedingly natural (universal graphs
omitting all cycles of odd order up to some fixed bound; or universal
graphs omitting a path and any further set of constraints). 

In this context, the pruning operation makes equally good sense, and
the corresponding monotonicity principle is valid. The decision problem makes more
sense in that context, and is equally open. But it is only at this
level of generality that encoding arguments make sense, so that one
can envision a ``soft'' proof of undecidability. At the same time,
the possibility of a complete ``list'' of favorable constraint sets
is viable. All known examples consist of a combination of some
very special constraints with a set closed under homomorphism. 
Whether this merely reflects our inexperience remains to be seen.
In any case, one can define a notion of ``critical set'' in general:
practically speaking, this would be a set which after pruning produces
a known example allowing a universal graph. In general, by determining
the critical sets which allow universal graphs, and are not in the
database of known examples, and iterating, one could arrive at
the correct answer in general, together with a proof of it.
It will be clear from our treatment here of an initial special case
arising in the case of a single constraint that this is a large task.
But the work of \cite{fk1}, for example, is encouraging. It follows
from the arguments given there that a finite set $\CC$ of $2$-connected
graphs allows a universal $\CC$-free graph if and only if the set is 
closed under homomorphism. Since any graph is built up in a reasonable
way from trees and $2$-connected graphs, we are off to a decent start.

It would also be interesting to pass on quickly 
to the case of a finite set of trees. This may be entirely reasonable.

\begin{center}\dots\dots\end{center}

The paper is organized as follows.  

In \S\ref{pruning}, we discuss the structural analysis of a general
finite connected graph $C$ as a ``tree of blocks'' (or 2-connected
components), a standard topic of graph theory which has a great deal
to do with the practical analysis of universality problems, and we
introduce the new idea which allows an inductive analysis of
universality problems according to the complexity of the underlying
tree.  Quite generally, universality problems can be reduced by this
method to canonical ``minimal'' cases, which we call ``critical.''
Here the constraint $C$ can be any finite connected graph, or in fact
any finite set of finite connected graphs.  
For our applications here, we will take $C$ to be a tree in later sections.

In \S\ref{monarchy}  we show that that the Tree Conjecture holds
for trees with a unique vertex of maximal degree. The method of
\S\ref{pruning} reduces this to the case treated in \cite{fk2}, and
this provides a nice illustration of the force of the reduction, as
well as disposing of a case that is best treated in isolation.
This serves as a template for more elaborate constructions.

In the following sections we prove the Tree Conjecture by making a very
coarse division of the critical cases into subcases according 
to the maximal vertex degree and the structure of the ``external''
vertices of maximal degree (those nearest the leaves).

We make little mention of algebraic closure in the remainder of the
paper, apart from an occasional observation. While the ability to
compute this operation is important when investigating a new example,
it would contribute little or nothing to the exposition. 
But the notion will nonetheless be quite
visible in all of our constructions, in
the form of infinite paths of very tightly linked
elements.
If we were only
interested in the issue of local finiteness we could shorten both our
constructions and our analysis considerably, and also reduce
the number of distinct cases considered. 

Note that all graphs dealt with here are finite or countable.
We may mention their countability for emphasis, on occasion.

There are many other universality problems involving infinite
forbidden subgraphs, infinite sets of finite forbidden subgraphs, 
and uncountable graphs, but none of our methods apply in such cases, except
possibly to the case of infinite sets of finite forbidden subgraphs,
(which includes such natural cases as graphs without circuits) 
where at the least new and mysterious phenomena arise, and the
decidability problem is ill-posed.


\section{Pruning trees, and other graphs \label{pruning}}

Our main objective in this section is to give a general inductive
method for treating universality problems involving a finite set
of finite connected constraints. 
It is based on the decomposition of a graph into {\it blocks}, or
{\it 2-connected components}, and the underlying tree structure that
results. First, we recall the definitions, which are standard.

Let $C$ be a graph, which more often than not will be taken to be
connected and nontrivial. We will assume in any case that $C$ 
contains no isolated vertices: every vertex lies on an edge. 
Define an equivalence relation on the
edge set $E(C)$ as follows. First, for $e,f\in E(C)$ write $e\sim f$ 
if $e$ and $f$ are either equal or lie on a cycle in $C$. Then extend
this relation to an equivalence relation $\approx$, the transitive
closure of $\sim$. A {\it block} of $C$ is the graph induced on the
set of vertices lying on the edges in a single equivalence class in
$E(C)$; this can consist of two vertices lying on a single edge. 
Blocks are 2-connected, that is they remain
connected after deletion of any vertex. 
Now a pair of blocks intersects in at most one vertex, and we
associate to the graph $C$ its {\it reduction} $\tilde C$ whose
vertices are the blocks of $C$, with two vertices connected just when
the blocks in question share a common vertex. The underlying structure
of  $\tilde C$ is a forest, and as we will be taking $C$ to be
connected, the reduction $\tilde C$ is even a tree. We call this 
{\it the underlying tree} of $C$. We believe this analysis is highly
relevant to our problem of universality. In fact, we believe the
following.

\begin{conjecture}[Solidity Conjecture]
If there is a $C$-free universal graph 
(in either the weak or strong sense), then the blocks of $C$
are complete.
\end{conjecture}

We call such a graph {\it solid}.
One could conjecture in general that the algebraic
closure operation should be unary. For the case of one constraint this
becomes the solidity conjecture, but for the case of multiple
constraints we do not know its precise content in graph theoretic terms. 
For example if $\CC$ is the class of all trees of order $n+1$ then
the algebraic closure operation is unary, and locally finite, since
each connected component of the graph has
order at most  $n$.

\begin{conjecture}[Reduction Conjecture]
If there is a $C$-free universal graph (in either of the two senses),
then there is a $\tilde C$-free universal graph, where $\tilde C$ is
the underlying tree of $C$.
\end{conjecture}

This is an instance of the ill-fated Monotonicity Conjecture discussed
in the introduction, which will be partially rehabilitated in the
present section. But it lacks any theoretical support, and is merely
plausible (and testable, fortunately).

Combining this with the Tree Conjecture, one gets a fairly precise
sense of what is expected, namely that beyond Komj\'ath's 2-bouquet,
similar bouquets, and 
some further substantial generalizations of that example, the class of
exceptional constraints allowing universal graphs (weakly or strongly)
should run out fairly soon.
In the background there is also the expectation, as noted earlier,
that a constraint allowing a weakly universal graph 
also allows a strongly universal one, 
though again not for any theoretical reason.

We move on from idle conjecture to something more
rigorous. 

\begin{definition}
Let $C$ be a connected graph consisting of more than
  one block.
\begin{enumerate} 
\item A pair $(B,u)$, where $B$ is a block 
and $u\in V(B)$, is called a {\em pointed block}.
\item A pointed block $(B,u)$ is an {\em attached leaf} of the graph $C$
if there is a block $B'$ of $C$ which represents a leaf in the
underlying tree of $C$, and a vertex $u'\in B'$ belonging to another block
of $C$, such that the pointed block $(B',u')$ is isomorphic to $(B,u)$.
\item A {\em minimal} attached leaf of $C$ is an attached leaf $(B,u)$
  such that there is no embedding of any other attached
  leaf $(B_1,u_1)$ into $(B,u)$ as a proper subgraph (that is, such an
  embedding must be an isomorphism).
\end{enumerate} 
\end{definition}

Observe that in the above, any block $B$ of $C$ which represents a
leaf of $\tilde C$ in fact meets exactly one other block of $C$ and
hence has a unique vertex of attachment.
Furthermore, any such pointed block containing a minimal number of
edges will be a minimal attached leaf, and similarly there are minimal
attached leaves among those with a minimal number of vertices.

What we wish to consider are the operations of {\it pruning} or {\it
  attaching} leaves of a particular minimal type, which we give in a
  slightly more general form.

\begin{definition} Let $C$ be a finite connected graph, 
$\CC$ a finite set of finite connected graphs, $G$ an
  arbitrary (in practice countable) graph, and $(B,u)$ a pointed block.
\begin{enumerate}
\item $C^-$ is the graph obtained from $C$ by {\em pruning} $(B,u)$:
  this means, for every attached leaf $(B',u')$ which can be embedded
  isomorphically into $(B,u)$, 
  we delete $V(B')\setminus\{u'\}$,
  and take the induced graph on the remaining vertices. Note that
  vertices lying in more than one block remain.
\item $G^\circ$ is the graph induced by $G$ on the set of those vertices
  $v$ of $G$ such that $G$ contains infinitely many copies of
  $(B,u)$, disjoint over $u$, with $u$ 
  identified with $v$.
\item $G^+$ is the graph obtained from $G$ by freely attaching
  infinitely many disjoint copies of $(B,u)$ to each vertex $v$ of
  $G$, with $u$ identified with $v$.
\item For sets $\CC$ of constraints, $\CC^-$ is the set of pruned
  graphs $C^-$ for $C\in \CC$. 
\end{enumerate}
\end{definition}

If greater precision is needed, we may write $C^-(B,u)$, $\CC^-(B,u)$,
$G^\circ(B,u)$, and  $G^+(B,u)$ instead. 

Now we come to the point. 

\begin{proposition}
Let $\CC$ be a finite  set of finite 
connected graphs and suppose there is a $\CC$-free
graph which is universal, either in the weak or strong sense. 
Let $(B,u)$ be an attached leaf of some graph $C$ in $\CC$,
and  $\CC^-=\CC^-(B,u)$ the result of pruning.
Then there is a universal $\CC^-$-free graph (in the same sense).
In fact, if $G$ is a universal $\CC$-free graph, then 
$G^\circ=G^\circ(B,u)$ is a universal $\CC^-$-free graph.
\end{proposition}

\begin{proof}
We have $\CC$, $(B,u)$, $G$, and $\CC^-,G^\circ$ 
as described, and we observe 
that $G^\circ$ is $\CC^-$-free. 

Now let $\Gamma$ be $\CC^-$-free and consider $\Gamma^+=\Gamma^+(B,u)$. 
Then
$$\hbox{$\Gamma^+$ is $\CC$-free}\leqno(1)$$

In any embedding of some $C\in  \CC$ 
into $\Gamma^+$ as a subgraph, $C$ will map
into a single connected component, and each block of $C$ will map into a
block of that component. 
The blocks of $C^-$ which do
not correspond to leaves of the associated tree 
$\tilde C$ go into $\Gamma$; and the
blocks of $C^-$ which do correspond to leaves of $\tilde C$ also go
into $\Gamma$, as none of them embeds into $(B,u)$ over the attaching
vertex. So an embedding of $C$ into $\Gamma^+$ would induce an
embedding of $C^-$ into $\Gamma$, and $(1)$ follows.

Now $\Gamma^+$ must embed in $G$, by $(1)$, either as a subgraph or as
an induced graph, as the case may be. Under such an embedding,
$\Gamma\includedin [{\Gamma^+}]^\circ$ will embed into $G^\circ$, 
either as a subgraph or as an induced graph, correspondingly. 
Our claim follows.
\end{proof}

In view of the importance of this result for our analysis, we make the
following definition in the case of a single constraint.

\begin{definition}
A finite connected graph $C$ is {\em critical} if the underlying tree
$\tilde C$ is neither a path nor a near-path, but for any 
type of attached leaf of $C$, the tree $\widetilde C^-$ associated with
the corresponding pruned graph is a path or near-path.
\end{definition}

\begin{corollary}
Suppose that $C$ is a finite connected graph whose underlying tree
$\tilde C$ is neither a path nor a near-path, and that there is a 
weakly or strongly universal $C$-free graph. Then there is an induced
subgraph $C'$ of $C$, which is critical, for which, correspondingly, a
weakly or strongly universal $C'$-free graph exists.
\end{corollary}

For the proof, one prunes $C$ repeatedly until it becomes critical. 

We have conjectured that there are no graphs with the properties of
the Corollary; and we see that it suffices to consider critical ones. 
We prove the Tree Conjecture in this framework by considering critical
trees. In this case attached leaves are essentially just leaves, or
rather edges connecting a leaf to its point of attachment, and pruning
amounts to the removal of the leaves (or to put it another way,
shortening all the external branches).



\section{Monarchy and Stardom\label{monarchy}}

In the present section we will prove the following.

\begin{theorem}\label{Monarchy:main}
Let $T$ be a tree with a unique vertex of maximal degree. Then the
following are equivalent.
\begin{enumerate}
\item There is a strongly universal $T$-free graph.
\item There is a weakly universal $T$-free graph.
\item The tree $T$ is a path or near-path.
\end{enumerate}
\end{theorem}

A further equivalence would be: algebraic closures of points
are either finite, or are two-way infinite paths without additional
edges. But as we have remarked, there is no need to bring in
the notion of algebraic closure explicitly in such cases.

The implication $(3\implies 1)$ requires argument, and is treated
in \cite{ct}. This has a completely different character from anything
we do here, lying on the positive side; all of our work here fills in
the gap on the negative side. We need to show $(\neg  3 \implies \neg
1)$. 

We distinguish two cases. Let $d$ be the maximal vertex degree in the
tree $T$. Then either $d\ge 4$ or $d=3$. 
One might expect this distinction to be significant, 
since the case 
$d=3$ includes the case of a near-tree, which at some point 
has to be singled out as an exception,
but there are other reasons for the case distinction as well.

In the present section, what will be important is the
behavior of a regular tree with vertex degree $d-1$, and more
precisely of its approximations, namely regular graphs of vertex
degree $d-1$  
and large girth. There is certainly a distinction to be observed here
between the case $d=3$ and $d\ge 4$. Later on, the issue will be
somewhat different. We will need to construct infinite graphs
from finite pieces while controlling the vertices of degree $d$, and
it is difficult to avoid introducing new vertices of degree $3$.

While this case distinction is not {\em always}
essential, it tends to play a role, and the case $d=3$ is the more
complicated of the two. On the other hand,
as we work with critical trees there is some
compensation in the form of improved control of the structure of the
tree in this case.

\subsection{Monarchs}

We begin with the generic case, $d\ge 4$,
and while dealing with this case we will encounter
all the issues that arise in any of the cases, as well as most of the
strategies for dealing with them.

\begin{proposition}
Let $T$ be a tree with a unique vertex of maximal degree $d$, with
$d\ge 4$. Then there is no weakly universal $T$-free graph.
\end{proposition}

\begin{proof}
Suppose toward a contradiction that $T$ is a counterexample of minimal
order. Then there is a weakly universal $T$-free graph, and
hence for the graph $T'$ derived from $T$ by pruning
(removal of its leaves) there is also a universal $T'$-free graph (\S\ref{pruning}). 
Now if $T'$ also has a vertex of degree $d$, then this violates the
choice of $T$ as a minimal counterexample. 
So the vertex $v$ of degree $d$ in $T$ must have at least one leaf as
a neighbor in $T$.

Now this turns out to be precisely the situation considered in
\cite{fk2}: a tree with a unique vertex $v$ of maximal degree $d$, which
is adjacent to a leaf $v'$, and with $d\ge 4$.
\end{proof}

Let us expand on this, as the same type of construction and analysis
is needed in general with a host of minor complications.
The following construction applies in the critical case of \cite{fk2}.

Let $\Gamma$ be a regular tree of degree $d-1$, or more generally a
regular graph of degree $d-1$ which is {\it tree-like} in the 
sense that the girth is large (larger than $2n$ with $n=|T|$).

As $\Gamma$ is regular of degree $d-1$ it is $T$-free (this part of
the argument blows up considerably as soon as we leave the domain of
monarchy), and we may vary the construction of
$\Gamma$, and in particular the cycle lengths occurring in $\Gamma$,
to give $2^{\aleph_0}$ graphs of this type. Here (and only here) we exploit
the hypothesis $d-1>2$.

The key property is the following.
\begin{texteqn}{(*)}
If $\Gamma$ is a subgraph of a $T$-free graph $G$, then $\Gamma$ is a
connected component of $G$, and is an induced subgraph.
\end{texteqn}

We will check this. Another way to phrase this claim
is that the vertices of
$\Gamma$ can acquire no new neighbors. 
If our aim is to refute
only {\em strong} universality, then
most of this argument drops out of the
picture, along with any preparations 
which may have been made for it.
There are no such preparations in the present case, but usually there
will be. 

It is immediate that $\Gamma$ is a connected component of $G$,
because as soon as one adjoins a new neighbor $u'$ to a vertex $u$ in
$\Gamma$, which is not already a vertex of $\Gamma$,
one gets an embedding of $T$ into the extended graph by
identifying $v$ with $u$, a leaf adjacent to $v$ with $u'$, and the
rest of $T$ with a suitable part of $\Gamma$, which locally (near $u$)
looks like a regular tree of degree $d-1$. 

Similarly, there can be no new edge between vertices $v,w$ whose
distance in $\Gamma$ is greater than $n=|T|$.
Local connections require more attention; this is also a
characteristic feature of the more complicated constructions later.

If we adjoin a new edge $(v,w)$ between two
nonadjacent vertices of $\Gamma$ which lie at distance at most $n$, 
then we may embed $T$ into the resulting graph as follows.
Let $P$ be the path linking $v$ to $w$ in $\Gamma$, and let $v''$ be the
neighbor of $v$ on $P$. Then $v''$ may play the role of $v'$, and the
part of $\Gamma$ remaining after deleting the component of
$\Gamma\setminus \{v,w\}$ containing $v''$, and with the edge $(v,w)$
adjoined, again will be regular of degree $d-1$, and girth 
greater than $n$, so the remainder of $T$ can be embedded over $v,v''$.

\includegraphics[height=1.6 in]{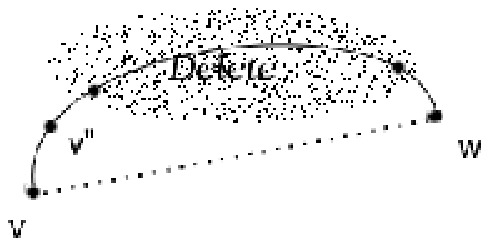}

Now it is impossible that all of these graphs $\Gamma$ could occur inside
a single countable graph among its connected components, so no countable
weakly universal $C$-free graph exists in this case.

For such constructions we require $T$-free graphs to which the
adjunction of a single edge will in many cases produce an embedding of
$T$, without having a detailed knowledge of the structure of $T$. 
In general the main mechanism for keeping the necessary control
involves paying attention to the distribution of vertices of degree
$d$ in $\Gamma$ (as there are none in this case, the distribution is
particularly transparent), and building up the minimal vertex degrees
to $d-1$. On the other hand the idea of making the graph closely resemble a
tree will have to be severely curtailed in general, 
and instead by considering
critical trees we will find we need much less control of $\Gamma$ to ensure
the necessary embeddings become available. Fortunately we do not need
to prevent the addition of arbitrary edges: it will be sufficient if we
can recover some invariants of $\Gamma$ once we know the restriction
of the embedding of $\Gamma$ to a suitable finite set.
When $\Gamma$ is a connected component of the ambient graph, we can
recover $\Gamma$ itself from the image of any vertex, but we need much
less than that.
         
\subsection{Stardom}

The method of the previous section clearly will not work with $d=3$.
So we now deal separately with this case, reducing to the
critical case and making use of a construction that looks closely at
the structure of the constraint tree. This is very reasonable, since
we still need to distinguish the exceptional cases.

We will refer to vertices of degree at least $3$ in a tree as {\it
  branch vertices}.

\begin{definition}
\par\indent
\begin{enumerate}
\item A {\em star} is a tree with a unique branch vertex, called its
  {\em center}.
\item For $n\ge 1$ and $d_1\ge \dots\ge d_n\ge 1$, the tree
$S(d_1,\dots,d_n)$ is the star formed by attaching paths of length
  $d_1,\dots,d_n$ to a central vertex.
\end{enumerate} 
\end{definition}

We will always take $n\ge 3$ here, to get a proper star. In this case the
star has a well-defined {\it center} and the maximal vertex degree is
$n$. Near-paths are stars $S(d_1,d_2,1)$.
Stars with $n\ge 4$ have been dealt with in the preceding subsection.

\begin{proposition} If $S=S(d_1,\dots,d_n)$ is a star and is not a
  near-path, then there is no weakly universal countable $S$-free graph.
\end{proposition}

\begin{proof}
Since the case $n\ge 4$ is covered by the previous Proposition, 
we will take $n=3$, and we will also take $S$ critical in the sense
of \S\ref{pruning}, which means
that we take
$$d_3=2,$$ 
so that pruning produces a near-path with the same center.

Let $H_0$ and $H_1$ be the following graphs. First, fix two vertices
$u_0,u_1$. To form $H_0$, adjoin two common neighbors $v_0,v_1$ to
$u_0$ and $u_1$, with $v_0$ and $v_1$ adjacent; this is $K_4$ with one
edge deleted. To form $H_1$, adjoin infinitely many common neighbors
$v_i$ to $u_0$ and $u_1$, and add an edge $(u_0,u_1)$, 
with no further adjacencies.

\includegraphics[width=4.5in]{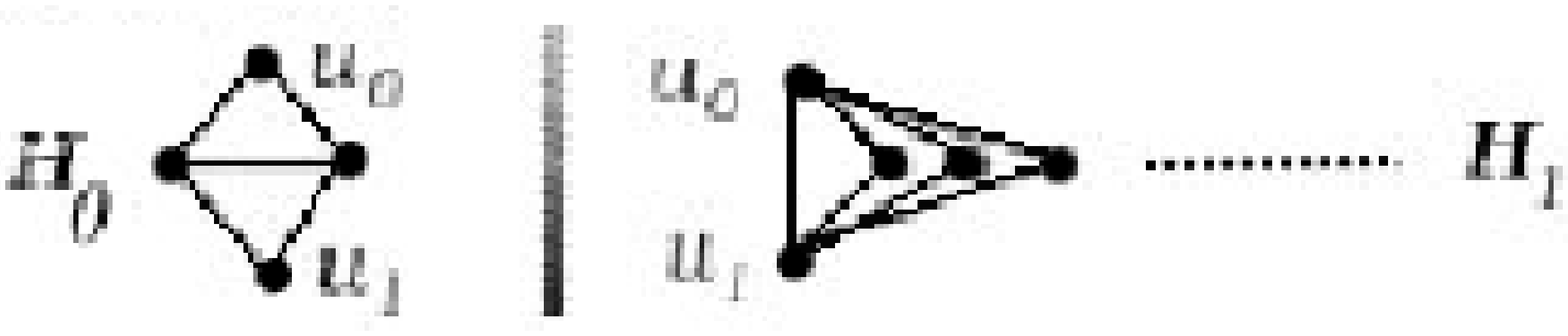}

Now for $\epsilon\in 2^{\Zz}$ a bit string, form a graph
$\Gamma^\epsilon$ as follows. Begin with an infinite independent set
$A$ of vertices $a_i$ ($i\in \Zz$). For each $i$, attach to the pair
$a_i,a_{i+1}$ a copy of $H_{\epsilon(i)}$ with $a_i$ and $a_{i+1}$
corresponding to $u_0$ and $u_1$.
Then $\Gamma^\epsilon$ is $S(2,2,2)$-free and in particular $S$-free.

Now we have to think about ``decoding'' $\Gamma^\epsilon$ when it is
embedded in a larger $S$-free graph as a subgraph.
As usual this involves getting some control over at least some of 
the additional edges adjoined in such an extension.

Define a relation
$R(u,v)$ on the vertices of any graph $G$ by the following condition: 
{\it $u,v$ lie in a copy of
$H_0$ or $H_1$, with $u$ and $v$ playing the roles of $u_0$ and $u_1$
respectively.}
By construction successive pairs $(a_i,a_{i+1})$ satisfy this relation
in any graph $G$ into which $\Gamma$ embeds. We must show that this
relation is not much affected by embedding into a larger $S$-free graph.
Our claim is as follows.

\begin{texteqn}{(1)}
If $\Gamma^\epsilon \includedin G$ and $G$ is $S$-free, then for each
vertex $a_i\in V(\Gamma^\epsilon)$, and each $v\in V(G)$, if
$R(a_i,v)$ holds then $v=a_{i\pm 1}$.
\end{texteqn}

Let $P$ be a path in $\Gamma^\epsilon$ containing all $a_i$, with
$d_P(a_i,a_{i+1})=2$ for all $i$. 

Now either $a_i$ and $v$ have infinitely many common neighbors,
or $a_i$ and $v$ play the roles of $u_0$ and $u_1$ in $H_0$ (or both).

If $v$ has infinitely many neighbors, then
$v$ must lie on the path $P$, as otherwise we may embed $S$ into 
the extension of $\Gamma$ by one of the new edges attached to $v$.
But there is some
freedom in the choice of the path $P$, and if $v$ cannot be pushed off
it by altering the path, then
$v$ must in fact be some $a_j$.
Now it is easy to see that if $|j-i|>1$ then there is an embedding of 
$S$ embeds into $G$, a contradiction.

So $v$ has finite degree in $G$ and thus
$a_i$ and $v$ must play the role of $u_0$ and $u_1$ in $H_0$. 
So they have common neighbors $w,w'$ which are adjacent. 

If $w\in A$, that is $w=a_j$ for some $j$, 
then $a_i$ and $a_j$ are adjacent in $G$, and easily $j=i\pm 1$. 
It follows easily that $w$ and $w'$ are
not both in $A$. So we may assume that $w$ is
not in $A$, and choose $P$ so that $w$ is not on $P$. Then 
$v$ must be on $P$. If $v=a_j$ for some $j$, then again by inspection
$j=i\pm1$ as claimed. So we may suppose that $v$ is not in $A$,
but lies on $P$, and that $P$ cannot be chosen to avoid both $v$ and
$w$. 
This means that 
$v$ and $w$ are the common neighbors of some pair
$(a_j,a_{j+1})$, and are the only such common neighbors, as otherwise
the path $P$ could still be moved.
But $w'$ is adjacent to $v$ and $w$, hence by
considering $a_j,a_{j+1}$ we see that $w'$ is also forced onto $P$,
and hence must be $a_k$ for some $k$. However, looked at from the
point of view of $a_k$, this is also impossible: $a_k$ becomes the
center of a copy of $S$.

So $(1)$ holds. We can now deduce the nonexistence of a weakly universal
$S$-free graph. Suppose toward a contradiction that $G$ is a weakly
universal $S$-free graph, and, of course, countable.  
For each $\epsilon$
choose an embedding $f_\epsilon$ of $\Gamma^\epsilon$ into $G$. Choose
a pair $\epsilon$, $\epsilon'$ for which these embeddings agree on the
successive vertices $a_0,a_1$ in $A$. It follows from $(1)$ that
the restriction of $f_\epsilon$ to A coincides with the restriction of
$f_{\epsilon'}$ to $A$.
But for some $i$, we
have $\epsilon(i)\ne \epsilon'(i)$, and thus the vertices $v_0,v_1$ in
$G$ which correspond to $a_i,a_{i+1}$ in both $\Gamma_\epsilon$ and
in $\Gamma_{\epsilon'}$ must occur on copies of both $H_0$ and $H_1$ in
$G$. This immediately provides an embedding of $S$ into $G$, and a
contradiction.
\end{proof}

The final argument is typical and will occur in some form in all
cases.  
As long
as the essential invariant $\epsilon$ of $\Gamma$ can be recovered
from the embedding (after fixing some points to get rid of shifts and
reflections along $A$), we can argue in this fashion.
On the other hand, if we are dealing with strong universality, 
there would be little to check at this point.
Still, even when dealing with induced subgraphs one has to check for
example that the relation $R(a_i,v)$ is not satisfied by new elements
of $G$, and indeed without this one would not even know
that the algebraic closure operation is nontrivial. So in this decoding
phase, the issues are similar regardless whether we deal with
local finiteness, strong universality, or weak universality, though
the degree of control needed to effect the decoding varies considerably.


\section{Toward the Tree Conjecture}

Our goal now is the following.

\begin{theorem}\label{TC1:main}
If $T$ is a finite tree with maximal vertex degree $d\ge 3$, 
and if $T$ has more than one vertex of degree $d$, then 
there is no weakly universal $T$-free graph.
\end{theorem}

In view of the result of the previous section, it suffices to 
prove this theorem in the critical case. So we record it in 
this form.

\begin{namedtheorem}{Theorem \ref{TC1:main}$'$}
If $T$ is a critical finite tree with maximal vertex degree $d\ge 3$, 
and if $T$ has more than one vertex of degree $d$, then 
there is no weakly universal $T$-free graph.
\end{namedtheorem}

Recall that in the critical case the pruned tree $T'$ is
either a path or a near-path. 

\subsection{A special case}

We first prove a considerably weaker result in which the basic construction
can be seen most simply, and without invoking the criticality hypothesis.

\begin{proposition}\label{tc1:stronguniversality}
Let $T$ be a tree with maximal vertex degree $d\ge 5$, and suppose
that every vertex of degree $d$ is adjacent to a leaf of $T$.
Then there is no strongly universal $T$-free graph.
\end{proposition}

We describe the construction of an uncountable family of countable 
$T$-free graphs,
and show that they cannot all be simultaneously embedded into a $T$-free
countable graph. There are three phases to this argument: $(a)$
Construction; $(b)$ $T$-freeness; $(c)$ Decoding
(i.e., 
analysis of the image of such a graph under embedding into a larger
$T$-free graph). The assumption that $d\ge 5$ simplifies
the construction, and the fact that we deal with strong universality
simplifies the decoding process by limiting the class of embeddings
considered. The arguments for $T$-freeness amount to saying that the
metric structure induced by $T$ on its vertices of degree $d$ does not
embed into the metric structure induced by our graphs on their
vertices of degree $d$ or more, and is typical of the analysis in
general.

The extra hypothesis on the vertices of degree $d$ is much stronger 
than what we actually require below, and much weaker than what one has
if one restricts attention to critical trees. Some form of this condition
is helpful in stage $(c)$. 

\begin{definition}
\par\indent
\begin{enumerate}
\item $V_1(T)$ is the set of vertices of degree $d$ in $T$, construed
  as a metric space with the induced metric. From this one can recover
 the tree structure induced on the convex hull of this set in $T$.
\item If $v\in T$, then a {\em $v$-component} of $T$ is a connected
  component of the graph resulting from deletion of $v$ in $T$.
\item A vertex $v$ in $V_1(T)$ is an {\em external} vertex of maximal
  degree if it is a leaf in the convex hull of $V_1(T)$ in $T$.
  Equivalently, at most one $v$-component of $T$ contains vertices of
  degree $d$.
\end{enumerate}
\end{definition}

For the proof of the Proposition we may assume that there are at least
two vertices in $T$ of maximal degree, as otherwise we apply
Theorem \ref{Monarchy:main}. All we really require for the proof of
this proposition is a single external vertex of degree $d$ adjacent to
a leaf.

\begin{construction}
Let $v_1$ be an external vertex of $T$ of maximal degree and let $C$
be the $v_1$-component of $T$ containing all other vertices of $T$ of
maximal degree. Let $H$ be the graph induced by $T$ on $C\union \{v_1\}$.

Let $P$ be a $(d-2)$-regular $2$-connected graph of very large girth.
For each vertex $u\in P$, attach a copy $H_u$ of $H$ to $u$ with 
$u$ corresponding to $v_1$. Call the resulting graph $\Gamma^P_0$.
For any vertex $u$ of degree less than $d-1$
in $\Gamma^P_0$, raise its degree to $d-1$ by adjoining
suitable trees (regular of degree $d-1$ except at the root, where the
degree is $d-1-\deg(u)$). 
Call the resulting graph $\Gamma^P$.
\end{construction}

\begin{lemma}
The graph $\Gamma^P$ is $T$-free.
\end{lemma}
\begin{proof}
Let $V_1(\Gamma^P)$ denote the set of vertices of $\Gamma^P$ of degree
at least $d$ construed as a metric space with the induced metric.
These vertices in fact have degee exactly $d$ and lie in the subgraphs
$H_u$, with $V_1(H_u)$ isometric to $V_1(T)\setminus \{v_1\}$.
It will suffice to show that there is no embedding of $V_1(T)$ into
$V_1(\Gamma^P)$ as metric spaces which is semicontractive in the
sense that distances do not increase. Note that this metric structure
is the same in $\Gamma^P_0$ and in $\Gamma^P$, so for the rest of this
argument one may as well think in terms of $\Gamma^P_0$.

Call a subspace $A$ of $V_1(T)$ {\em isolated} if it satisfies the
following condition:
\begin{texteqn}{Iso}
For every subspace $A_0$ of $V_1(T)$ isometric with $A$ and every
embedding $f$ of $V_1(T)$ into $V_1(\Gamma^P)$, the image $f(A_0)$
is contained in some single $V_1(H_u)$. 
\end{texteqn}
A better term might be ``indecomposable'' but we wish to emphasize
here that the subgraph $P$ is avoided, something which will be less
clear in subsequent constructions.

Now if there are no such embeddings $f$ then $V_1(T)$ itself is
isolated, but if there are any such embeddings then $V_1(T)$ is not
isolated. So let us assume there are such embeddings 
and let $A$ be an isolated subspace of $V_1(T)$ of maximal order.
By our assumption $A\ne V_1(T)$, and
as points are isolated, $A$ is nonempty.

Choose a pair $(B,v)$ with $B$ a subspace of $V_1(T)$ isometric with
$A$, and with $v\in V_1(T)\setminus B$, and furthermore with
$$\delta=d(f(v),f(B))$$
minimized over all such pairs, and all semicontractive embeddings $f$
of $V_1(T)$ into $V_1(\Gamma^P)$.
By the maximality of $A$, $B'=B\union
\{v\}$ is not isolated, and thus there is an embedding
$f:V_1(T)\to V_1(\Gamma^P)$ such that the image $f(B')$ meets at least
two distinct sets of the form $V_1(H_{u})$ with $u\in P$. 
But $f(B)$ is contained in one such set $V_1(H_u)$, and thus $f(v)$ is
contained in another. 
Now let $\tilde B$ be the subspace of $V_1(T)$ 
corresponding to $f(B)$ under the identification of $H$ with $H_u$,
and observe that $\tilde B$ is isometric with $A$, 
that $v_1\notin \tilde B$, and 
that $d(f(v_1),f(\tilde B))\le d(v_1,\tilde
B)=d(u,f(B))<d(f(v),f(B))$, the latter point 
in view of the structure of the metric on $\Gamma^P$.
So $d(f(v_1),f(\tilde B))<\delta$, contradicting our choice of
$\delta$ as minimal.
\end{proof}

In later arguments we will use some of this metric terminology while
formulating the main argument directly in terms of graph embeddings. 
We simply wished to emphasize here that the obstruction really is
captured by the metric structure on the vertices of high degree. But
the argument in general will depend a little more at the graph structure,
particularly near vertices of $P$.

The second question to take up is a kind of rigidity (or decoding) for 
$\Gamma^P$ when considered as an induced subgraph of a 
general $T$-free graph.

\begin{lemma}
Let $G$ be a $T$-free graph containing $\Gamma^P$ as an induced
subgraph.
Then for any vertex $u\in P$, the neighbors of $u$ in $G$ are its
neighbors in $\Gamma^P$.
\end{lemma}
\begin{proof}
Here we recall the assumption that the girth of $P$ is large, and thus
the local structure of $P$ near the vertex $u$ is exactly that of 
a $(d-2)$-regular tree. We also use the assumption that $v_1$ has a
neighbor which is a leaf.

What we need to show is the following: the graph $\Gamma^P_{u,v}$ obtained
by adjoining one new neighbor $v$ to $u$ contains a copy of the tree
$T$. One begins the construction of a suitable embedding by
taking the map identifying $H$ and $H_u$, in which $v_1$ corresponds
to $u$. Now some leaf $v_1'$ adjacent to $v_1$ may correspond to $v$.
It remains to embed the remaining $d-2$ $v_1$-components of $T$ into
$\Gamma^P$, making use of $P$ and some of the trees attached at the
end of the construction.

Now each of the remaining $(d-2)$ $v_1$-components of $\Gamma^P$
contains a tree closely resembling a 
$(d-1)$-regular tree (except for its root, adjacent to
$v_1$, whose degree in the component is one less); any cycles will be
the cycles of large girth allowed in $P$. So it is easy to see that
there is no obstruction to the completion of our embedding.
\end{proof}

Our two lemmas prove the proposition. This is a general principle; let
us check it in this case.

\begin{proofof}{Proposition \ref{tc1:stronguniversality}}
Suppose toward a contradiction 
that $G$ is a countable strongly universal $T$-free graph.
Consider embeddings $f_P:\Gamma^P\to G$. We have uncountably many
isomorphism types of $P$ available, as we may control the cycle
lengths that appear. So there must be at least two nonisomorphic
graphs $P,Q$ whose images $f_P(P)$ and $f_Q(Q)$ meet in a vertex $u$.

Now consider the subgraph $G_0$ of $G$ induced on the vertices of degree
exactly $d-1$ in $G$. This contains $f(P)$ and $f(Q)$, by our lemma.
In particular the connected component of $u$ in $G_0$ contains
$f(P)$ and $f(Q)$. Now observe that the connected component $G_u$ of $u$ in
$G_0$ is contained in $f(\Gamma^P)$. Otherwise, we would have an edge
$v,v'\in G$, with $v$ in $f(\Gamma^P)$ and $v'\notin f(\Gamma^P)$, and
with $v$ of degree $(d-1)$ in $G$. But $v$ is already of degree at least $d-1$
in $\Gamma^P$ and thus $v'\in f(\Gamma^P)$, a contradiction.

So $G_u$ is contained in both $f(\Gamma^P)$  and $f(\Gamma^Q)$. 
Now the only nontrivial block in $f(\Gamma^P)$ is $f(P)$ and similarly
for $Q$, so $f(P)=f(Q)$ and as these are induced subgraphs of $G$, 
the graphs $P$ and $Q$ must be isomorphic, a contradiction.
\end{proofof}

To put the matter briefly, our second lemma proves that $P$ is
recoverable from finite data, and any countable graph contains only
countably many candidates for such data, so if we have uncountably
many candidates for $P$ then can be no universal graph (strongly or
weakly, depending on the strength of the recoverability lemma).

We can see that this argument is going to require signficant
adaptation as we remove the simplifying hypotheses. For
$d=4$ the graph $P$ becomes a path (this is why we call it $P$,
actually) and for $d=3$ we will again use a path $P$, but we will 
have to look considerably farther into the structure of $T$ to
find a suitable way to extend $P$ without creating a copy of $T$.
The most extreme case was treated earlier: the case of stars. There
the construction looks very little like the one just given, though
the path $P$ is still visible.

In the decoding phase most of the weight was borne above by
the fact that we dealt only with strong universality. So in most
cases we will have to modify the construction to ``block'' the
adjunction of at least some potential new edges to our graphs
$\Gamma$. 

The method used in the proof of this proposition
lies at the core of most of
our subsequent constructions and proofs.

Criticality will be important to ensure an adequate supply of leaves,
and it also simplifies the embedding argument made in the decoding
phase. In the critical case the embedding argument in the proof of our
decoding lemma would work in $\Gamma^P_0$ as well as in $\Gamma^P$,
but on the other hand the full decoding argument was only given in the
proof of the proposition and this argument actually needs $\Gamma^P$
rather than $\Gamma^P_0$, so the saving is not very great here.

In some delicate cases criticality may also help in checking
that our graphs $\Gamma$ are $T$-free. Most of our constructions
place additional vertices of degree $d$ on the graph corresponding to
$P$ here, so the analysis must become more precise.

\subsection{Amalgamation and the parameter $\ell$}

Before entering into a detailed consideration of how the foregoing
construction may be adapted to deal with the question of weak
universality, we may consider some general points that are relevant to
the decoding process and give some indication of what additional
structural features of the constraint tree $T$ are relevant in
general. This leads to a certain proliferation of cases, handled by a
unified method but with considerable variation from case to case.

One such parameter, and an important one, is the maximal degree $d$.
The case $d=4$ turns out not to be much more troublesome, in comparison
with $d\ge 5$. While $P$ is just a 2-way infinite path in this case it
turns out that there are some variations available in the ``attachment''
procedure that passes from $P$ to $\Gamma^P$ and we can again arrive
at uncountably many variations on each theme. Ultimately the same will
apply when $d=3$ but not so simply.

But there is a second parameter which comes into play in the decoding
phase. When our graph $\Gamma$ embeds into a larger $T$-free graph $G$ it
may acquire new edges between its own vertices, and this ``noise''
threatens to make recovery of $\Gamma$ from $G$ impossible.
However, what will be true is that the vertices of $P$ will remain of
finite degree, and that just as we considered the vertices of degree
$d-1$ in the previous subsection, consideration of the vertices of
finite degree is generally useful. 

The critical observation is the
following: if $\hat T$ is the result of amalgamating two copies, or
even infinitely many copies, freely over the vertex set $V_1(T)$ then
$V_1(\hat T)=V_1(T)$ as a metric space (but here $V_1(\hat T)$ means
all vertices of degree at least $d$, not exactly $d$).

Thinking back to the ``attachment graph'' $H$ of the previous
subsection, if $\hat H$ is the corresponding subgraph of $\hat T$,
then this suggests the idea of using $\hat H$ in place of $H$ and
going on as before. This would be sound apart from one fatal flaw: 
all the vertices along $P$ are likely to acquire infinite degree in
the process. 
This makes it extremely likely that $T$ will embed in the graph so
constructed, and also makes the recovery of $P$ highly improbable.
In short, everything needed is destroyed.

There is one case in which this flaw is not actually present: if the
external vertex $v_1$ of degree $d$ is adjacent to a vertex $v_0$ of
degree $d$ then $v_1$ will have no new neighbors in $\hat H$.
Really what we are doing in this case is working
with the subgraph $H_0$ of $H$ obtained by deleting $v_1$ and the
corresponding amalgam $\hat H_0$, then attaching $\hat H_0$ to $v_1$
by an edge.

In general, we must consider the parameter
$$\ell=d(v_1,v_0)$$
where $v_0$ is the closest vertex of degree $d$ to $v_1$.
We want to treat the part of $H$ based at $v_0$ as the ``attachment graph, ''
and take further pains to deal with the path from $v_1$ to $v_0$,
which we think of as potentially running along $P$.
It turns out that the relevant case division is as follows: 
$\ell\ge 3$; $\ell=2$; $\ell=1$ with the first two cases similar and
the last of a different character.

Again, once one enters into this kind of more precise construction,
the structure of the tree $T$ between $v_0$ and $v_1$ plays a major
role; as the pruned tree $T'$ will be a path or nearpath one hopes
that the corresponding part of $T'$ will be just a path, though
a few exceptional configurations must be treated separately.

So our case division comes out something like the following.

\begin{enumerate}
\item[I] $\ell\ge 2$:
\begin{enumerate}
\item [A] $\ell \ge 3$; \quad B \  $\ell =2$.
\end{enumerate}
\item [II]$\ell=1$:
\begin{enumerate}
\item [A] $d\ge 4$; \quad B \  $d=3$ (with various subcases).
\end{enumerate}
\item [III] Left-over nearpaths
\begin{enumerate}
\item[A] $d\ge 4$; \quad B \  $d=3$.
\end{enumerate}
\end{enumerate}

We will ultimately list the cases differently for reasons of
convenience, but the logical structure is properly reflected above. 
The full list of cases actually used is recapitulated at the end.


\section{Case I: $\ell\ge 2$}

We take up the proof of Theorem \ref{TC1:main}$'$.
We deal with an external vertex
$v_1$ of maximal degree, and a closest vertex $v_0$ of maximal degree,
with $\ell=d(v_0,v_1)$.
 
\subsection{Case IA: $d\ge 4$, $\ell\ge 3$}

\begin{texteqn}[0.70]{Case IA}
$T$ has a vertex $v_0$ of maximal degree $d\ge 4$ such that some
  $v_0$-component $C$ of $T$ contains a unique vertex $v_1$ of degree
  $d$, and  $\ell=d(v_0,v_1)\ge 3$. 
Either $T'$ is a path, or else $T'$ is a near-path
whose center does not lie in the $v_0$-component $C$. 
\end{texteqn}

We allow $v_0$ to be the center of $T'$. 

\begin{construction}
Let $H=T\setminus C$.
Let $P$ be a $(d-1)$ regular graph of large girth, and let
$\Gamma^P$ be the result of adjoining a vertex $b$ adjacent to all
vertices of $P$ and then attaching a copy $H_b$ of $H$ with $b$
corresponding to $v_0$
\end{construction}

\begin{lemma}
$\Gamma^P$ is $T$-free.
\end{lemma}
\begin{proof}
Since the number of vertices of degree $d$ in $T$ is greater than the
number of vertices in $H_b$ which have degree at least $d$ in
$\Gamma^P$, any embedding $f$ of $T$ into $\Gamma^P$ has to carry at least
one vertex of degree $d$ into $P$.

If exactly one vertex $u$ of degree $d$ in $T$ corresponds to a vertex
of of $P$, then $b$ must also correspond to a vertex of degree $d$ in
$T$ under the embedding $f$. Now the vertex $f(u)\in P$ has degree
exactly $d$ in $\Gamma^P$ and hence $b$ must be one of the neighbors
of $u$ in $f(T)$. As $d(u,b)=1$ it follows that the diameter of the
set of vertices of degree $d$ in $f(T)$ is less than its diameter in
$T$, so this is not an isomorphism.

Thus there are at least two vertices $u,v$ of degree $d$ in $T$ whose
images under $f$ lie on $P$. Again, the vertex $b$ must occur as a
neighbor of $u$ and $v$ in the image $f(T)$ and therefore $d(u,v)=2$.
Since $\ell\ge 3$ there must be other vertices of degree $d$ in $T$
(or it would be enough for our purposes to assume this, if $\ell=2$).
It follows easily that $T$ must be a nearpath with center
corresponding to $b$. Then it is easy to see that the diameter of 
the set of vertices of degree $d$ in $T$ is greater than the diameter
of the corresponding set in $f(T)$, a contradiction.
\end{proof}

\begin{lemma}
Suppose that $G$ is a $T$-free graph containing $\Gamma^P$
and that $u$ is a vertex of $P$. Then any neighbor $v$ of $u$
in $G$ is either a neighbor of $u$ on $P$, or a vertex of $H_b$.
\end{lemma}
\begin{proof}
In view of the structure of $T'$, $v_1$ is adjacent to a leaf of $T$.
If $v\notin \Gamma^P$ then we look for an embedding of $T$ into $G$
in which $u$ represents $v_1$, $v$ represents such a leaf, and $H_b$
represents $H$ with $b$ corresponding to $v_0$. The path from $v_0$ to
$v_1$ can run along $P$. As $P$ has locally the structure of a
$(d-1)$-regular tree, the extension to the remainder of $T$ is possible.

If $v\in \Gamma^P$ then we may suppose $v\in P$. Suppose that $v$ is
far from $u$ in the metric on $P$. Then we proceed as in the case when
$v$ is not in $\Gamma^P$. Now suppose $v$ is close to $u$, and the
girth of $P$ is large relative to the distance $d(u,v)$.
Then there is a unique shortest path $L$ from $v$ to $u$. Let $u'$ be
the neighbor of $u$ on $L$. Consider the graph obtained from $P$ by
deleting the rest of $L$ (between $u'$ and $v$) and adjoining the edge
$(u,v)$. If $P$ had been a $(d-1)$-regular tree then this new graph
would also be a $d-1$-regular tree, but in fact it is a
$(d-1)$-regular graph of high girth.
In any case, using $u'$ to represent a leaf adjacent to $u$ we may
again proceed as in the first case to embed $T$ into $G$.
\end{proof}

\begin{proposition}
In Case IA there is no (countable) universal $T$-free graph.
\end{proposition}
\begin{proof}
Otherwise we find ourselves considering embeddings
$f:\Gamma^P\to G$, $g:\Gamma^Q\to G$ with $P$ and $Q$ nonisomorphic,
and with the images $f[H_b^P]$ and $g[H_b^Q]$ identical, 
and also $f(P)$ meets $g(Q)$.
Now we look at the graph $G_0$ obtained by deleting $f[H_{b^P}]$.
In this graph, $f(P)$ and $g(Q)$ are connected components,
and the induced structure from $G_0$ is the original structure on $P$
or $Q$. As $f(P)$ meets $g(Q)$, the images coincide and the graphs are
isomorphic. 
\end{proof}

\subsection{Case IB: $d=3$, $\ell\ge 3$}

\begin{texteqn}[0.70]{Case IB}
$T$ has a vertex $v_0$ of maximal degree $d=3$ such that some
  $v_0$-component $C$ of $T$ contains a unique vertex $v_1$ of degree
  $d$, and  $\ell=d(v_0,v_1)\ge 3$. 
Either $T'$ is a path, or else $T'$ is a near-path
whose center is not $v_1$.
\end{texteqn}

\begin{construction}
Take an infinite path $P$ and partition it into successive
(alternating) finite intervals $P_i$, $Q_i$ for $i\in \Zz$ of lengths
$p_i,q_i$ respectively,
satisfying the following conditions.
\begin{enumerate}
\item $p_i\ge 3\ell-3$
\item $q_i=\ell$
\end{enumerate}

Let $C$ be the $v_0$-component of $T$ containing $v_1$ and
$H^0=T\setminus C$. Let $H$ be the graph obtained by amalgamating
two copies of $H^0$ freely over the vertices of $H^0$ of degree $d$ in
$T$ (this includes $v_0$), and adjoining additional vertices to bring
up the degree of any vertex in $V_1(H)$ to $\infty$.

Adjoin vertices $b_i$ adjacent to all vertices in
$P_i$, for all $i$, and attach a copy $H_i$ of $H$ to $b_i$ with
$b_i$ corresponding to $v_0$.
Call the result $\Gamma=\Gamma^\epsilon$ where $\epsilon$ is the sequence
$(p_i)_{i\in \Zz}$.
\end{construction}

\begin{lemma}\label{IB:Tfree}
$\Gamma$ is $T$-free.
\end{lemma}\
\begin{proof}
Let $X_i$ be the set of vertices in $P_i\union H_i$ of degree at least
$d$ in $\Gamma$, and let $Y_i$ be the set of vertices of $H_i$ of
degree at least $d$ in $\Gamma$. 

Call a subspace $A$ of $V_1(T)$ {\em indecomposable} if under every
embedding of $T$ into $\Gamma$, the image of any subspace of
$V_1(T)$ isometric with $A$ lies in one set $X_i$; call $A$ {\em
isolated} if every such image lies in one set $Y_i$.

We will write $V_1$ for $V_1(T)$ throughout.

\begin{texteqn}{(1)}
There is a nonempty isolated metric subspace of $V_1$.
\end{texteqn}

Consider a geodesic path $P^*=(p_0,\dots,p_n)$ in 
$V_1$ of maximal length subject to the condition $d(p_i,p_{i+1})\le
\ell$ for $i<n$, and among all such maximize
$d(p_0,p_1)+d(p_{n-1},p_n)$.

We claim that this geodesic path is isolated. First,
as the $H_i$ are widely separated ($q_i=\ell$) this path is
indecomposable. 

Now consider any $A\includedin V_1$ isometric to $P^*$ and embedded
into $X_i$ by an embedding $f$ of $T$ into $\Gamma$. We claim that
$f[A]$ lies in $Y_i$. If not, writing $A=(a_0,\dots,a_n)$ we have 
$u=f(a_i)\in P_i$ for some $i$.

Now $u$ has degree $d$ in $\Gamma$ and hence the neighbors of 
$a_i$ in $T$ map onto the neighbors of $u$ in $\Gamma$, including
$b_i$. So any vertices of $A$ which map into $P_i$ share a common
neighbor in $T$, and hence lie at distance $2$ in $T$. As $A$ is a
geodesic path, and corresponds to points on a path in $T$, 
there are at most two such points in $A$. 

If there are two such points in $A$ then as $b_i$ lies between their
images, the other vertices of $A$ map into $P_i$ rather than
$H_i$. But as this is impossible, we find that $|A|=2$ in this case,
and as $d(a_0,a_1)$ is maximized, that 
$$\ell=2$$ 
as well, contradicting our current assumptions.

So a unique point of $A$ maps into $P_i$. Then the other points of $A$ map
into $H_i$ and are linked to $u$ by a path through $b_i$. 
In particular $u$ must correspond to an endpoint of $A$. We may
suppose $u=f(a_n)$. Now the path $(f(a_0),\dots,f(a_{n-1}))$
corresponds to a path $(a_0',\dots,a_{n-1}')$ 
in $T\setminus C$, and $d(a_{n-1}',v_0)=d(a_{n-1},b_i)<
d(a_{n-1},u)\le \ell$. Hence if $a_{n-1}'\ne v_0$, 
this path may be lengthened to a path
$(a_0',\dots,a_{n-1}',v_0,v_1)$ satisfying our conditions and
contradicting the maximality of $n$. We conclude that
$a_{n-1}'=v_0$. But then the path
$(a_0',\dots,a_{n-1}'=v_0,v_1)$ satisfies our conditions with an
increase in the distance between the last two vertices, again
contradicting maximality. 
This last contradiction completes the proof of $(1)$.

Now we consider a subspace $A\includedin V_1$ which is isolated, and
maximal. One possibility is that $A=V_1$, but as $|V_1|>|V_1(H)|$,
this would mean that there are no embeddings of $T$ into $\Gamma$, as
we claim. So for the remainder of the argument we suppose 
that there are such embeddings, and we 
aim at a contradiction.

We claim that there is some subspace $A'$ of $V_1$ isometric to $A$,
and some embedding $f$ of $T$ into $\Gamma$ such that for some $i$ 
we have
$$b_i\in f[A']\leqno(2)$$

We choose a pair $A'\includedin V_1$ and $v\in V_1\setminus A'$ so
that $A'$ is isometric to $A$ and $d(v,A')$ is minimized, and we
consider the space $B=A'\union \{v\}$ inside $V_1$. By our choices of
$A$, this is not isolated. Take $B'=A''\union \{v'\}$ isometric
to $B$ inside $V_1$ so that there is an embedding $f$ of $T$ into
$\Gamma$ for which $f[B']$ does not go into any $H_i$.
As $A'$ is isolated, $f[A'']$ goes into some $H_i$, and by hypothesis
$f(v')\notin H_i$, so 
$d(b_i,A'')<d(v',A'')$. By the choice of $v'$, we must have $b_i\in
f[A'']$ and thus $(2)$ is achieved.

Now we repeat the general thrust of the first part of the argument. 
We consider a geodesic path $P^*=(a_0,\dots,a_n)$ 
which can be attached to $A'$ at the
point $u_0$ corresponding to $b_i$ in $f[A']$, so that with the natural
metric the extension $A^*=A'\oplus_{u_0} P^*$ is isometric with a subspace
of $V_1$, and so that $d(a_i,a_{i+1})\le \ell$ for all $i<n$, and we
first maximize $n$, then maximize $d(a_{n-1},a_n)$.
By condition $(2)$ one possibility is to take $P^*=(v_0,v_1)$ with
$v_0$ corresponding to some point $u_0$ in $A'$, so $|A^*|>|A'|$.

By the maximality of $A'$,  the space $A^*$ cannot be isolated.
It is certainly indecomposable since $A'$ is indecomposable, in view
of the metric structure of $A^*$ and $\Gamma$. 
So there is an embedding $f$ of $T$ into $\Gamma$ taking a copy of 
$A^*$ (which we continue to call $A^*$) into 
$X_i$ for some $i$, but not into $H_i$. Here $A'$ goes into $H_i$
and the geodesic path $P^*$ does not, though the endpoint $a_0$ does. 
Arguing as in the first instance we see that one end of $P^*$ goes
into $P_i$ and the rest of $A^*$ goes into $H_i$. 
But then we adjust $P^*$ as before and obtain a final contradiction:
if $b_i$ is not in the image of $P^*$ we lengthen the path $P^*$,
while if $b_i$ is in the image of $P^*$ we move $a_n$ farther away. 

Retracing our steps from this contradiction, we see that in fact the
maximal isolated space $A$ must be $V_1$, and thus $\Gamma$ is
$T$-free. 
\end{proof}

We can now enter the decoding phase.

\begin{lemma}
If $\Gamma$ is contained in the $T$-free graph $G$, 
then for any vertex $u$ of $P$ we have the following.
\begin{enumerate}
\item The degree of $u$ in $G$ is finite.
\item Any neighbor of $u$ of finite degree in $G$ lies on $P$, and is
  one of the neighbors of $u$ in $\Gamma$.
\item If $i$ is chosen so that $d(u,P_i)$ is minimal, then any
  neighbor of $u$ in $G$ off $P$ lies in $H_i$.
\end{enumerate}
\end{lemma}
\begin{proof}
As $d=3$ the $v_0$-component containing $v_1$ consists of a path
with one edge adjoined at $v_1$.

Suppose $(u,v)$ is an edge with $u\in P$, $v\in G$ and $v\notin
\Gamma$. Then we embed $T$ into $G$ with $u$ corresponding to $v_1$ 
and $v$ corresponding to a leaf of $T$ adjacent to $v_1$. It suffices
to notice that if $i$ is chosen to minimize $d(u,P_i)$ then there is a
path of length $\ell$ from $u$ to $b_i$, and $b_i$ can play the role
of $v_0$.
The same applies if $v\in H_j$ with $j\ne i$.
It remains to consider vertices of finite degree in $H_i$ and vertices
on $P$.

Now as $H_i$ is obtained by free amalgamation of two copies of $H^0$
over the vertices of degree $d$, and the vertices originally of degree
$d$ are transformed into vertices of infinite degree,
the vertices of finite degree in $H_i$ may be treated just like
vertices outside $\Gamma$.  

Suppose therefore that $v\in P$. Then we must consider the possibility
that $v$ lies along our intended path $L$ from $u$ to $b_i$. 
If $u\notin P_i$ then we may substitute for $L$ a path beginning with 
$(u,v)$. The main point to consider is the possibility $u,v\in P_i$,
but $u,v$ are nonadjacent.

Then after deleting $u,v$ from $P_i$ there remain  $p_i-2\ge
(3\ell-5)>3(\ell-2)$ vertices, and at least one of the three resulting
subintervals in $P_i$ contains at least $\ell-1$ vertices. Furthermore
as the roles of $u$ and $v$ are now symmetric, we may suppose that $u$
is an endpoint of such an interval.  
As we have a path of length $\ell$ from $u$ to $b_i$ which does not
pass through $v$, we may take $u,b_i$ to correspond to $v_1,v_0$
respectively, and use $v$
to represent a leaf adjacent to $v_1$. One may also use 
an additional vertex of $P_i$ to represent a leaf adjacent to $b_i$ 
as 
$$p_i\ge 3\ell-3>\ell+1$$
(For $\ell=2$ one would just add this inequality as a
restriction on $p_i$, but there are other difficulties in that case.) 

This proves the lemma in all cases.
\end{proof}

\begin{proposition}
In Case IB there is no (countable) universal $T$-free graph.
\end{proposition}
\begin{proof}
If $G$ is a countable universal $T$-free graph we
can find embeddings $f:\Gamma\to G$ and $f':\Gamma'\to G$
with $\Gamma=\Gamma^\epsilon$, $\Gamma'=\Gamma^{\epsilon'}$ 
nonisomorphic and with $f(P)$ meeting $f'(P')$, where $P'$ is the copy
of $P$ associated with $\Gamma'$.

We look at the graph $G_0$ induced by $G$ on its vertices of finite
degree. This contains $f(P)$ and $f'(P')$. Furthermore the connected
component of their intersection coincides with both, as an induced
subgraph. So the path $f(P)=f'(P')$ is an induced subgraph of $G$.
Call this path $P^G$. 

Now we consider the following relation $R(a,b)$ in $G$:
$$\hbox{$a,b$ are vertices of $P^G$ with a common neighbor in $G$
  which is not in $P_G$}$$
Writing $a=f(a_0)$ and $b=f(b_0)$ with $a_0,b_0\in P$, choose $i$ and
$j$ to minimize $d(a_0,P_i)$ and $d(b_0,P_j)$ respectively.
By our lemma, if  $R(a,b)$ holds in $G$ then $i=j$ (and in particular
  $i,j$ are uniquely determined). 

Now consider the equivalence relation generated by the relation $R$ on 
$P^G$. If $A$ is an equivalence class with representative $a=f(a_0)$,
and $i$ is chosen to minimize $d(a_0,P_i)$, then $A$ is contained in
the image of the set $\{v:d(v,P_i)\le \ell/2\}$. Furthermore $A$
either contains the image of $P_i$ or is disjoint from it. If we
consider equivalence classes of order at least $3\ell-3$, 
these will contain the corresponding set $f(P_i)$ as otherwise the
size of $A$ would be bounded by $\ell<3\ell-3$. 
So we can now identify the equivalence classes containing the
$f(P_i)$. These must also be the equivalence classes containing the
$f'(P_i')$ and thus one finds $|p_i-p_i'|\le \ell$, up to a shift or
and possible reflection of indices. By restricting the
allowed sequences $(p_i)$ somewhat one may ensure that this forces
$p_i=p_i'$ for all $i$, and thus a contradiction.
\end{proof}

\subsection{Case IC: $\ell= 2$}

\begin{texteqn}[0.70]{Case IC}
$T$ has a vertex $v_0$ of maximal degree $d\ge 3$ such that some
  $v_0$-component $C$ of $T$ contains a unique vertex $v_1$ of degree
  $d$, and  $\ell=d(v_0,v_1)= 2$. 
Either $T'$ is a path, or else $T'$ is a near-path
whose center is not $v_1$.
\end{texteqn}

\begin{construction}
Let $H^0$ be $T\setminus C$ and let $H^1$ be the amalgam of two copies
of $H^0$ over 
the set $V_1(H^0)$ of vertices in $H^0$ corresponding to vertices of
degree $d$ in $T$. 
Extend $H^1$ freely so as to raise the degree of vertices of degree at
least $d$ to $\infty$, and to raise all vertex degrees to at least
$d-1$. Call the result $H$.

Let $P$ be either a $(d-1)$-regular graph of large girth, if $d>3$,
or else a two-way infinite path, if $d=3$, and suppose that there is a
family of paths $P_i$ contained in $P$, of order $p_i$, satisfying
\begin{enumerate}
\item $p_i=2$ or $3$ all $i$;
\item Every vertex of $P$ is either on some $P_i$ or adjacent to one
  of its vertices;
\item No two vertices on distinct paths $P_i,P_j$ are adjacent.
\end{enumerate}

Attach to each interval $P_i$ a vertex $b_i$ adjacent to its vertices,
and attach a copy $H_i$ of $H$ to $b_i$ with $b_i$ playing the role of
$v_0$. 

Call the resulting graph $\Gamma$.
\end{construction}

\begin{remark}
One has uncountably many possibilities for the structure of $P$
if $d\ge 4$, and for the sequence $p_i$ if $d=3$.
\end{remark}

Let us verify this in case $d\ge 4$. In this case first choose a path
$P^*$ and intervals $P_i$ on $P$ meeting our conditions. Then add
edges to $P^*$ whose endpoints lie outside all the $P_i$, keeping the
girth high, and raising the vertex degrees to $d-1$.

\begin{lemma} $\Gamma$ is $T$-free.
\end{lemma}

\begin{proof}
This is the usual metric argument. We will sketch the main points.

We begin by considering a maximal geodesic path $A=(a_0,\dots,a_n)$
embedding into $V_1=V_1(T)$ with successive distances at most $\ell$,
and suitably maximized.  We have $n\ge 1$ by the case assumption.

We claim that this path is isolated relative to $\Gamma$
in our usual sense.

Note that as $p_i\le 3$ and the vertices of $P_i$ have degree $d$, 
with a common neighbor $b_i$, 
no embedding of $T$ into $\Gamma$ can carry
two vertices of degree $d$ into the same interval $P_i$. 
So as in the proof of Lemma \ref{IB:Tfree} it follows that
the path $A$ is isolated, and
after that the argument is relatively formal.
One considers a maximal isolated subspace $B$ of $V_1(T)$ which may be
supposed proper, and one finds 
that there must be some embedding in which some $b_i$ is in the image of 
an isometric copy of $B$, after which one can 
attach another such geodesic path to $B$
and arrive at a contradiction;
the gap between distinct intervals $P_i$ becomes relevant again when
the geodesic path argument is repeated at the end.
\end{proof}

\begin{lemma}
Let $\Gamma$ be embedded in the $T$-free graph $G$.
Then any vertex $u\in P$ is of finite degree in $G$, 
and its neighbors in $G$ of finite degree are exactly its neighbors in
$P$. 
\end{lemma}
\begin{proof}
Suppose first that $(u,v)$ is an edge of $G$ and $v$ does not occur
in $\Gamma$. Then $u$ will play the role of $v_1$ in $T$, with
$v$ an adjacent leaf. One easily finds a path of length $2$
connecting $u$ to some $b_i$, which will play the role of $v_0$, and 
one extends this to an embedding of $T$ into $\Gamma$, with $H_i$
absorbing $T\setminus C$ while the component $C$ itself, apart from
one leaf attached to $v_1$, embeds into $P$.

The same construction applies whenever the vertex $v$ is not needed to
complete the embedding of $T$ into $\Gamma$, and in particular only
finitely many vertices $v$ require attention, so the vertices of $P$
certainly continue to have finite degree. 

We now need to consider only the case in which $v$ is in $\Gamma$ and
has finite degree (in $G$, and in particular in $\Gamma$).
If $v\in H^1$ then $v\notin V_1(H^0)$ and hence by the amalgamation
process used to construct $H^1$ such a choice of $v$ 
cannot block anything. 
It is also possible that $v$ lies in $H_i$ but off $H^1$, 
but this is essentially the same situation; indeed, we could have
extended $H^0$ before amalgamating, and then $H$ would be just the
result of the final amalgam!

So all that really concerns us is the possibility that $v$ is on $P$.
But as we have seen previously, we 
can delete most of the path from $u$ to $v$ along $P$,
just retaining the neighbor $v'$ of $u$ along that path, and then $v'$
acts as a ``new'' vertex with respect to the revised version of $P$.
\end{proof}

\begin{proposition}
In Case IC, there is no weakly universal $T$-free graph.
\end{proposition}
\begin{proof}
We need to show that we can recover information, either  about 
the sequence $(p_i)$ or the structure of $P$ (if $d\ge 4$), from
an embedding of $\Gamma$ into a larger $T$-free graph $G$, given
the image of a vertex in $P$. For any $u\in P$, the neighbors of $u$
of finite degree in $G$ are its neighbors in $P$, and thus at least 
the graph $P$ can be recovered from $G$.
If $d\ge 4$ there is sufficient flexibility in the structure of $P$
to complete the argument. 

Suppose therefore that $d=3$, and the induced structure on $P$ is an
ordinary two-way infinite path.
We must decode some information about the numbers $p_i$.
Consider the graph $P^*$ on the path $P$ whose edges are the
edges of $P$ whose endpoints have a common neighbor in $G$.
Then vertices adjacent in $\Gamma$ to
distinct intervals $P_i,P_j$ are not adjacent in $P^*$, 
and each nontrivial
connected component of $P^*$ consists of an interval $P_i$ with
possibly one or both of its neighbors on $P$ adjoined. 

If one looks at a 
long interval $L$ in $P$, one can use $P^*$ to 
count accurately the number of intervals
$P_i$ which meet $L$, and the number of vertices involved, and find the
average value of $p_i$ over the interval. This is sufficient to
discriminate between substantially different parameter sequences,
taken to be constant over long intervals.
\end{proof}

\section{\label{sharp} Case II: $\ell=1$, $d\ge 4$}

In this case we have an external vertex of maximal degree adjacent
to another vertex of maximal degree.

\subsection{$\ell=1$, $d\ge 5$}

\begin{texteqn}[0.70]{Case IIA}
$T$ has a vertex $v_0$ of maximal degree $d\ge 5$ such that some
$v_0$-component $C$ of $T$ contains a unique vertex 
$v_1$ of degree $d$, and $v_0,v_1$ are adjacent.
Either $T'$ is a path, or $T'$ is a near-path whose center 
does not lie in the $v_0$-component $C$.
\end{texteqn}

This case is essentially the same as the illustrative example treated
in Proposition \ref{tc1:stronguniversality}.

\begin{construction}
Let $H^0=T\setminus C$ and let $H^1$ be the amalgam of two copies of $H^0$ 
with over its vertices of degree $d$ in $T$.
Adjoin vertices adjacent to the vertices of
degree at least $d$ in $H^1$ in order to make their degrees infinite.
This yields an attachment graph $H$.

Take a two-way infinite path $P$ and attach a copy $H_i$ of
$H$ to a neighbor $b_i$ of each vertex $a_i\in P$,
with $b_i$ playing the role of $v_0$.

Finally, bring up the vertex degrees along $P$ to exactly
$d-1$ (initially these are of degree $3$, and $d>4$). Do this by
adding additional edges to $P$, but keep the girth of the graph
induced on $P$ very large.

The result is called $\Gamma$.
\end{construction}

\begin{lemma}
$\Gamma$ is $T$-free.
\end{lemma}
\begin{proof}
As there are no vertices of degree $d$ on $P$, the notions of
indecomposability and isolation coincide in this case.

Take a maximal indecomposable subspace $A$ of $V_1$.
If $A\ne V_1$, 
take a pair
$A',v$ with $A'$ isometric to $A$ and contained in $V_1$, $v\in
V_1\setminus A'$, and $d(v,A')$ minimized. Let $A_1=A'\union \{v\}$
and as $A_1$ is decomposable take an isometric copy
$A_1'=A''\union \{v'\}$ of $A_1$ in $V_1$ and an 
embedding $f:T\to \Gamma$ which
witnesses this. Then $f[A'']$ will be contained in some $H_i$ and
$f(v')$ will lie in a different $H_j$, and farther than $b_i$. 
By the minimality of $d(v',A'')$ we have $b_i$ in the image of $A''$. 
Pulling this back into $T$, we have an isometric copy $A^*$ of 
$A$ in $T\setminus C$ containing $v_0$. So $A^*\union \{v_1\}$ is also
a subspace of $V_1$, and as $A^*$ is indecomposable and $v_1$ is
adjacent to a vertex of $A^*$, $A^*\union \{v_1\}$ is also
indecomposable. This however contradicts the maximality of $A$.

Thus $A=V_1$ is indecomposable. However $|V_1(T)|=|V_1(H)|+1$ and thus
$V_1(T)$ cannot embed in a copy of $H$. So there are no such embeddings,
and $\Gamma$ is $T$-free.
\end{proof}

\begin{lemma}
Let $G$ be a $T$-free graph containing $\Gamma$, and $u\in P$. Then
any neighbor $v$ of $u$ of finite degree in $G$ is on $P$, and is a
neighbor of $u$ in $\Gamma$.
\end{lemma}
\begin{proof}
Once the vertex $u$ acquires degree $d$, we extend to the neighboring
copy of $H_i$ and a path along $P$, together with suitable neighbors
(all distinct by our restriction on the girth). 

If $v$ lies in the copy $H_i$ of $H$ associated with $u$, 
and has finite degree, then it is a vertex duplicated in the
construction of $H^1$ (or one of the additional neighboring vertices
added at the end, which present no problems). 
Such a vertex cannot block the embedding of $T$.

There remains the possibility that the 
vertex $v$ lies on $P$ and is not a neighbor of 
$u$ in $P$. Then as in the proof of Proposition
\ref{tc1:stronguniversality} we use the neighbor of $u$ on the path
toward $v$ to represent a leaf adjacent to $v_1$, and use the
additional edge $(u,v)$ to replace $P$ by a similar $(d-1)$-regular
graph of large girth.
\end{proof}

\begin{proposition}
In Case IIA there is no weakly universal $T$-free graph.
\end{proposition}
\begin{proof}
Given a $T$-free graph extending $\Gamma$ and the image of a point in
$P$ we recover the set $P$ and the graph induced on it by $G$, which
is the same as the graph induced on $P$ by $\Gamma$. As $d>4$ there is
some latitude in the structure of this graph (in particular, in the
lengths of circuits in the graph) and thus we can recover uncountably
many different invariants.
\end{proof}

\subsection{$\ell=1$, $d=4$}

\begin{texteqn}[0.70]{Case IIB}
$T$ has a vertex $v_0$ of maximal degree $d=4$ such that some
$v_0$-component $C$ of $T$ contains a unique vertex 
$v_1$ of degree $d$, and $v_0,v_1$ are adjacent.
Either $T'$ is a path, or $T'$ is a near-path whose center 
does not lie in the $v_0$-component $C$.
\end{texteqn}

\begin{construction}
We vary the preceding construction. With the same attachment graph
$H$, we take an infinite path $P$ and divide it into consecutive
intervals $P_i$ of length $p_i=1$ or $2$, where furthermore $p_i=1$
with rare, and widely spaced, exceptions.

The vertex $b_i$ is attached to the interval $P_i$ and the graph $H_i$
is appended to it. There is no final decoration phase since $d=4$ and all
vertices of $P$ have degree $d-1=3$ already. The only variability is
in the sequence $(p_i)$.
\end{construction}

This is $T$-free as before; any sequence $(p_i)$ with $p_i=1$ or $2$
would be suitable at this stage, as the duplication of the neighbor of
$b_i$ on $P$ has no substantial effect.

\begin{lemma}
If $G$ is a $T$-free graph containing $\Gamma$ and $u\in P$, then
$u$ has finite degree in $G$, and every neighbor $v$ of $u$ of finite
degree in $G$ is on $P$.
\end{lemma}
\begin{proof}
The ``widely spaced'' condition on the $p_i$ is used here.

The argument goes as before except that one should pay some attention to
vertices $u$ lying on or near an interval $P_i$ of length $2$, as the neighbor
$u'$ of $u$ in $P_i$ may be unsuitable for our purposes, having degree
$d-1$ but sharing a neighbor with $u$. 
As the $v_0$-component $C$ may
continue on past $v_1$, and the next vertex after $v_1$ may have
degree $d-1$, this constrains us to working on $P$ on one definite
side of $u$. 

If the vertex $v$ lies on
$P$ on this preferred side of $u$, but beyond the neighbor of $u$ according to
$\Gamma$, it potentially blocks the embeddding of $T$. But
then we can use the neighbor of $u$ on that side to represent a leaf of $T$
and take a path through $u,v$ and along $P$ to complete the
construction. 
\end{proof}

\begin{proposition}
In Case IIB there is no weakly universal $T$-free graph.
\end{proposition}
\begin{proof}
The path $P$ is recoverable and has no extra structure. 
Adjacent points belonging to distinct intervals $P_i$ can have no
common neighbors, as there are then two attachment graphs $H_i$ available
and the common neighbor could lie in at most one of them.
Therefore the intervals $P_i$ and the numbers $p_i$ are also visible
in any $T$-free graph containing $\Gamma=\Gamma^\epsilon$, and the
usual argument applies.
\end{proof}



\section{\label{d=3} Case III: $\ell=1$, $d=3$}

Since we now deal with the case $d=3$ all branch vertices of $T$ have
maximal degree. These cases require a slightly finer consideration of 
the structure of $T$, bringing in the location of a third branch
vertex, assuming there is one.

\subsection{A special case}

\begin{texteqn}[0.70]{Case IIIA}
$T$ contains exactly two branch vertices $v_0$ and  
$v_1$, which are adjacent, and of degree $3$.
\end{texteqn}

\begin{construction}
Take a two-way infinite path $P$ and divide it into intervals
$P_i,Q_i$ 
which are alternately of length $p_i= 3$ and $q_i=1$ or $2$.
Adjoin a common neighbor $b_i$ to each interval $P_i$.
Call the result $\Gamma$.
\end{construction}

\begin{proposition}
In case IIIA, there is no weakly universal $T$-free graph.
\end{proposition}
\begin{proof}
The graph $\Gamma$ 
is $T$-free, and under any embedding into a larger $T$-free
graph $G$, $\Gamma$ will be a connected component of $G$.
Hence no countable $T$-free graph contains all possible variants 
of $\Gamma$.
\end{proof}

This is important, because we need to move somewhat
further away from near-paths before we can make 
suitable constructions of any generality.

\subsection{Three adjacent branch vertices}

Here is one case which is sufficiently far from the near-path 
case to be handled uniformly.

\begin{texteqn}[0.70]{Case IIIB}
The maximal vertex degree is $3$. $T$ contains a sequence of three
adjacent branch vertices $v_1,v_0,v_2$ with $v_1$ external and
adjacent to a leaf. Some $v_0$-component of $T$ is a path attached to $v_0$. 
\end{texteqn}

\begin{construction}
Call the $v_0$-components of $T$ $C,C_1,C_2$ where $v_i\in C_i$ for
$i=1,2$. By hypothesis $C$ is a path. 
Let $H$ be the graph obtained from $C_2$
by freely amalgamating infinitely many
copies of $C_2$ over the subset $V_1(C_2)$ consisting of its branch 
vertices in $T$. 

Take an infinite path $P$ partitioned into intervals $P_i,Q_i$ of lengths
$p_i=2$ or $3$ and $q_i=1$ and adjoin a vertex $b_i$ adjacent to the
vertices of $P_i$.
Attach a copy $H_i$ of $H$ to $b_i$ with $b_i$ playing the role of
$v_2$. This yields $\Gamma$.
\end{construction}

\begin{lemma}
$\Gamma$ is $T$-free.
\end{lemma}
\begin{proof}
Let $P^*$ be a path of maximal length consisting of adjacent branch
vertices of $T$. We claim that $P^*$ is isolated with respect to
embeddings of $T$ into $\Gamma$.  Certainly $P^*$ is indecomposable, 
and for any of embedding of $T$ into $\Gamma$ which takes a vertex
$a$ of $P^*$ into $P_i$, as these vertices have degree $3$ and a
common neighbor, only one of them can be in $f[P^*]$. Hence $a$ must
be an endpoint of $P^*$ and its neighbor $b$ in $P^*$ must correspond
to $b_i$. In particular if $P_0$ is the path with $a$ deleted, then
$f[P_0]$ is a path in $H_i$ terminating at $b_i$, which corresponds to
a path of adjacent branch vertices in $C_2$ terminating at $v_2$. 
Such a path can be extended by $v_0,v_1$ and contradicts the maximality
of the length of $P^*$. So $P^*$ is isolated.

Now let $A$ be a maximal isolated subspace of $V_1=V_1(T)$. 
Assuming that there is in fact some embedding of $T$ into $\Gamma$,
then $A\ne V_1$, and then on formal grounds as we have seen 
in earlier arguments, there is an embedding of $T$ into $\Gamma$ which
carries an isometric copy of $A$, which we will continue to call $A$,
into some $H_i$ with $b_i$ included in the image. 
But then looking at this inside $T$ it gives an isometric copy of $A$,
say $A'$, containing $v_2$ but not $v_0$. So consider the longest path
$\tilde P$ consisting of adjacent vertices which can be
attached to the metric space $A$ at the corresponding vertex $v$, subject
to the restriction that the extended space $A\oplus_v \tilde P$ with
its natural metric embeds into $V_1$. This is visibly indecomposable
and easily seen to be isolated by the same sort of analysis with which
we began. This then contradicts the maximality of $A$ and completes
the analysis.
\end{proof}

\begin{lemma} 
For any $T$-free graph $G$ containing $\Gamma$, and any vertex $u\in P$, 
the degree of $u$ is finite in $G$, and 
the neighbors of $u$ of finite degree in $G$ and in $\Gamma$ coincide.
If $u\in Q_i$ for some $i$ then its neighbors in $G$ are on $P$.
\end{lemma}
\begin{proof}
Let us first see how to embed $T$ in $G$ if $u$ has a neighbor $v$ not
in $\Gamma$. We take $u'$ adjacent to $u$ and lying in one of the
$P_i$. We use the sequence $b_i,u',u$ to represent the sequence
$v_2,v_0,v_1$. The graph $H_i$ disposes of any need to think
about the component $C_2$. There is room for the path $C$ on the far
side of $u'$ along $P$. The vertex $v$ reperesents a leaf adjacent to
$v_1$, and the rest of $T$ consists of a path attached to $v_1$, which
can lie along $P$.

From this it follows that the vertices of $P$ have finite degree in
$G$.  Now suppose $v$ is a neighbor of $u$ of finite degree in $G$,
and in particular $v$ is a vertex of $\Gamma$. As usual if $v\in H_i$
then this does nothing. So we may suppose $v\in P$, and $v$ is
nonadjacent to $u$.

If $v$ and $u'$ lie on opposite sides of $u$ then we use the neighbor
of $u$ on the side of $v$ to represent a leaf adjacent to $v_1$, and
use the continuation of the path $(u,v)$ along $P$ to complete the
embedding with no further interference.

If $u$ is on $Q_i$ for some $i$ then there are two choices for $u'$ so
we can fall directly into the previous case, and the analysis applies to
any neighbor $v$ of $u$ in this case.

So we may suppose $u$ is on $P_i$.
If $v$ and $u'$ are on the same side of $u$, the only obstruction
arises if $v$ is adjacent to $u'$. If $v\in Q_{i\pm 1}$ we can
interchange $u$ and $v$, so there remains only the case in which
$v,u',u$ are the three points of some $P_i$. In this case we may take
$b_i,v,u$ to represent the sequence $v_2,v_0,v_1$ and use $u'$ as a
neighbor of $u$.
\end{proof}

\begin{proposition} 
In Case IIIB there is no weakly universal $T$-free graph.
\end{proposition}
\begin{proof}
It follows at once from the preceding lemma that in any 
$T$-free graph $G$ containing one of our graphs $\Gamma$,
we can recover $P$ as well as enough information about the neighbors
of $P$ to determine the sequence $(p_i)$ up to reflection and translation
from an element of $P$. So the customary argument applies.
\end{proof}

Under the assumption that $T'$ is a path we have seen that
we may suppose that external branch vertices are adjacent to
branch vertices. With the last two cases out of the way
there must be at least four branch vertices, with 
each outer pair adjacent.

\subsection{Two adjacent branch vertices, $\ell\ge 3$}

\begin{texteqn}[0.70]{Case IIIC}
The maximal vertex degree is $3$. 
$T$ contains a sequence of three successive
branch vertices $v_0,v_1,v_1'$ with $v_1'$ external and
adjacent to a leaf, $v_1$ adjacent to $v_1'$,
and $\ell=d(v_0,v_1)\ge 3$, 
where $v_0$ is the closest branch vertex to
$v_1$ other than $v_1'$. 
Either the pruned tree $T'$ is a path, 
or a near-path
with the center not in the $v_0$-component containing $v_1$. 
All external branch vertices of
$T$ are adjacent to branch vertices. 
\end{texteqn}

\begin{construction}
Let $v_2$ be the vertex lying between $v_1$ and $v_0$ at
distance $2$ from $v_1$.
Let $C$ be the $v_2$-component of $T$ containing $v_1$.
Let $H$ be the result of amalgamating infinitely many copies of
$T\setminus C$ over the set consisting of 
its branch vertices together with 
the path from $v_2$ to $v_0$, extended to give each vertex
strictly between $v_2$ and $v_0$ infinite degree.

Take a path $P$ broken into intervals $P_i,Q_i$ of lengths
$p_i=3$ and $q_i=1$ or $2$.
Adjoin a vertex $b_i$ adjacent to the vertices of $P_i$,
adjoin a  vertex $c$ adjacent to all $b_i$, and 
attach $H$ to $c$ with $c$ playing the role of $v_2$.
Call the result $\Gamma$.
\end{construction}

\begin{lemma}
$\Gamma$ is $T$-free.
\end{lemma}
\begin{proof}
We look first at how an adjacent pair of branch vertices
can be embedded into $\Gamma$. Consider a subtree of $\Gamma$
consisting of two adjacent branch vertices $v,v'$
of degree three with
their neighbors.

Suppose $v\in P$. Then $v\in P_i$ for some $i$ and by inspection
$v'=b_i$, with $c$ belonging to the subtree as one of the neighbors of 
$v'$. 

Now consider an embedding of $T$ into $\Gamma$ and more particularly 
the images of a pair of branch vertices consisting of
an external branch vertex and its neighboring
branch vertex. If all such images miss $P$, then
the diameter of the convex hull in $T$ of the branch vertices
is larger than the diameter of the graph into which they
can embed.

Similarly if one of these pairs of adjacent branch vertices maps
to $vb_i$ with $v\in P_i$ and the other pair maps into $H$, the diameter
is still slightly too small as the distance from $v_0$ (in $H$)
to $b_i$ is $\ell-1$. 

Finally if both pairs of adjacent branch vertices correspond to
pairs of the form $vb_i$ and $v'b_j$ then $c$ occurs 
in the embedding as a neighbor of both $b_i$ and $b_j$
and the whole tree has only four branch vertices, with
the interior pair lying at distance $\ell=2$, which
contradicts our case hypothesis.
\end{proof}

\begin{lemma} 
For any $T$-free graph $G$ containing $\Gamma$, and $u\in P$,
the neighbors of $u$ of finite degree in $G$ are its neighbors in 
$P$, and possibly a vertex $b_i$ if $u$ is either in $P_i$, or 
else in $Q_i$ or $Q_{i-1}$ and adjacent to a vertex of $P_i$,
with $q_i=2$ or $q_{i-1}=2$ respectively.
\end{lemma}
\begin{proof}
Bearing in mind that any pair of vertices on $P$ with
one in $P_i$ and the other adjacent to it 
are candidates for the role of $v_1$ and $v_1'$ respectively in 
an embedding of $T$, we see first that we cannot adjoin
any new vertices as neighbors of $u$, secondly that the vertices 
in $H$ which are of finite degree are not available to serve
as neighbors as they were duplicated in the amalgamation
process, and thirdly that for $q_i=1$ as there are two
vertices adjacent to $u$ and lying in $P_i$ or $P_{i+1}$ respectively, 
there can be no new neighbors of $u$ in that case.
Of course if $u\in Q_i$ or $Q_{i-1}$ with
$q_i=2$ our statement allows for $b_i$ as a new neighbor 
and the only other vertex which would be a plausible candidate
for a new neighbor of $u$ would be the next vertex 
beyond the immediate neighbor of $u$, but in that case
this new neighbor of $u$ in $P_i$ could serve as an alternate
candidate to play the role of $v_1$.

After all this there remains the possibility that $u\in P_i$
and that $u$ has a new neighbor along $P$, not already
adjacent to it in $\Gamma$---notably, $u$ and this new neighbor
could be endpoints of $P_i$. However here one uses the new neighbor
of $u$ in the role of $v_1'$ and the old neighbor of $u$ in $P_i$
represents a leaf in the embedding.
\end{proof}

\begin{proposition} 
In Case IIIC there is no weakly universal $T$-free graph.
\end{proposition}

\begin{proof}
We have freedom in the choice of the $q_i$ and so
we need only check that we have possibilities for decoding.

Given an embedding of $\Gamma$ into a $T$-free graph $G$
we look at the graph $G_0$ induced on  vertices of finite degree
in $G$ and then we look at the graph $G_1$ induced on vertices
of degree at most three in $G_0$.

The vertices of $P$ occur on a path in $G_1$.
Their neighbors in $G_1$ consist of $P$ and those vertices
$b_i$ which are of finite degree in $G$ and have
no neighbors in $G$ of finite degree other than those in $P_i$. 
In particular the connected component of $G_1$ containing
$P$ is a subgraph of $\Gamma$. One cannot necessarily recover
the path $P$ itself since the midpoint of $P_i$ and $b_i$ have
similar properties, but we claim that we may recover the sequence $q_i$,
which is sufficient.

The 2-connected blocks of $G_1$ consist of certain edges of $P$
together with the induced subgraph on $P_i\union \{b_i\}$ whenever
$b_i\in V(G_1)$. In the latter case the endpoints of $P_i$ can be
recovered from the $2$-block. Let $Q$ be the subgraph of $P$ obtained
by deleting those vertices occurring as the midpoint of an interval
$P_i$ for which $b_i\in V(G_1)$. Then the graph $Q$ can be recovered
from $G_1$, with those pairs consisting of endpoints of some interval
$P_i$ distinguished. From this one can recover the sequence $(p_i)$
(up to shift and reflection).
\end{proof}

\subsection{Two adjacent branch vertices, $\ell= 3$}

\begin{texteqn}[0.70]{Case IIID}
The maximal vertex degree is $3$. 
$T$ contains a sequence of three successive
branch vertices $v_0,v_1,v_1'$ with $v_1'$ external and
adjacent to a leaf, $v_1$ adjacent to $v_1'$,
and $\ell=d(v_0,v_1)=2$, 
where $v_0$ is the closest branch vertex to
$v_1$ other than $v_1'$. 
Either the pruned tree $T'$ is a path, 
and both external branch vertices of
$T$ are adjacent to branch verticesm, or it is a near-path
with the center not in the $v_0$-component containing $v_1$.  
\end{texteqn}

In most cases the construction for the previous case
will work. With $\ell=2$, the vertex $c$ in the previous
construction becomes identified with $v_0$.

But in the proof that $\Gamma$ is $T$-free, we may encounter
an exception precisely in the case when the external 
vertices correspond to vertices of $P_i,P_j$ for some $i,j$
and their neighbors correspond to $b_i,b_j$. 
In this case $c=v_0$ occurs as a common neighbor of both
and the structure of the tree is determined: 
it has exactly four branch vertices separated by one vertex of degree
$2$. The only uncertain point concerns the precise lengths
of the paths emanating from the external path vertices.

Since the previous construction definitely fails in this
case we use a variant. 

\begin{construction}
Suppose $T$ is as described (four branch vertices,
adjacent in pairs, separated by one vertex of degree $2$).
Take a two-way infinite path $P$ divided into intervals
$P_i,Q_i$ of lengths $p_i=2$ and $q_i=1$ respectively.
Attach vertices $b_i$ adjacent to the vertices of $P_i$.

Attach a path $(b_i,c_i,c_i')$ to each $b_i$ and attach
an infinite family of infinite rays to each of the vertices
$c_i,c_i'$
(really what interests us is to have two infinite
rays at $c_i'$ and to ensure that the vertices $c_i$,
$c_i'$ have infinite degree).

Now take a maximal subset $S$ of $\Zz$ containing no adjacent
pairs in $\Zz$; in other words, if $S$ is arranged as
a sequence $(n_i:i\in \Zz)$ then $n_{i+1}-n_i$ is $2$ or $3$ for
all $i$. Note that there are many such sets.
Give each vertex $b_i$ ($i\in S$) infinite degree.

Call the result $\Gamma$.
\end{construction}

\begin{lemma}
$\Gamma$ is $T$-free.
\end{lemma}
\begin{proof}
Consider subtrees of $\Gamma$ consisting of two adjacent branch
vertices and their neighbors.
Such a tree either has its branch vertices off $P$ entirely,
or has branch vertices of the form $v,b_i$ with
$v\in P_i$ and $i\in S$.

As the distance (in $\Zz$)
between distinct elements of $S$ is greater than $1$, the
distance between such pairs in $\Gamma$ is greater than
$2$, unless they lie in the part of $\Gamma$ attached to
a single $P_i$. But here the diameter is too small.
\end{proof}

\begin{lemma}
For any $T$-free graph $G$ containing $\Gamma$, and $u\in P$, 
the neighbors of $u$ in $G$ of finite degree are its neighbors
in $P$ together with $b_i$ if $u\in P_i$ and $i\notin S$.
\end{lemma}
\begin{proof}
Observe that $b_i$ has finite degree if $i\notin S$ 
since for some adjacent $j=i\pm 1$ we have $j\in S$
and it follows easily that giving $b_i$ infinite degree produces an 
embedding of $T$. 

The rest is clear by inspection; adjoining a new vertex as
a neighbor of $u$
produces an embedding of $T$ into $\Gamma$ directly,
and as usual 
there are no serious candidates of finite degree in $\Gamma$.
\end{proof}

\begin{lemma}
For any $T$-free graph $G$ containing $\Gamma$, and $u=b_i$, 
with $i\notin S$,
the neighbors of $u$ in $G$ of finite degree are its neighbors
in $P$.
\end{lemma}
\begin{proof}
Taking $j=i\pm 1$ in $S$, so that $b_j$ has infinite degree,
then as noted in the previous argument
adjoining a new neighbor of $b_i$ would give an
embedding of $T$ directly involving $b_i$ and $b_j$.

Therefore the only new neighbors which $b_i$ might acquire
are those lying along the image of this embedding.
Now the embedding passes in part along rays which
have been duplicated and since these rays have been duplicated
none of their elements can be a neighbor of $b_i$.
The remaining elements in the image of the embedding
either have infinite degree or lie on $P$ and therefore
are either irrelevant or excluded already in the previous lemma.
\end{proof}

\begin{proposition} 
In Case IIID there is no weakly universal $T$-free graph.
\end{proposition}
\begin{proof}
The path $P$ with its neighbors $b_i$ ($i\in S$) can
be recovered from any $T$-free graph containing $\Gamma$
together with one vertex of $P$, so the set $S$ can
be determined up to a shift (or exactly if 
two specific vertices of $P$ are fixed).

As usual there are uncountably many possibilities for $\Gamma$
and only countably many realized in any particular countable graph. 
\end{proof}


\section{\label{nearpath} The Tree Conjecture}

\subsection{Taking stock}

We review the analysis from the very beginning.
A minimal counterexample $T$ to the Tree Conjecture
will have the following properties.
\begin{enumerate}
\item $T$ is neither a path nor a near-path.
\item The pruned tree $T'$ is a path or a near-path
\end{enumerate}

Let $d$ be the maximal vertex degree. Then $d\ge 3$.
If there is a unique vertex of degree $d$ then
Theorem \ref{Monarchy:main} applies. So we
suppose the contrary.

\begin{enumerate}
\item [3.] There are at least two vertices in $T$ 
of degree $d$.
\end{enumerate}

Suppose first
\begin{enumerate}
\item [A]\quad $T'$ is a path.
\end{enumerate}

Let $v_0,v_1$  be a pair of vertices of degree $d$ with
$v_1$ external and with $v_0$ the closest vertex of degree $d$ to $v_1$.
Let $\ell=d(v_0,v_1)$.

If $\ell\ge 2$ then one of Cases IA-C applies.
Hence we suppose
\begin{enumerate}
\item[A1.] Any external vertex of maximal degree is adjacent to a
  vertex of maximal degree. 
\end{enumerate}

Then if $d\ge 4$ one of Cases IIA-B applies.
So we suppose
\begin{enumerate}
\item[A2.] $d=3$
\end{enumerate}

Then the four cases IIIA to IIID 
cover the remaining possibilities, 
case IIIA when there are exactly two branch vertices 
and one of Cases IIIB, IIIC, IIID otherwise.

Now suppose
\begin{enumerate}
\item [B] \quad $T'$ is a near-path with center $v_*$.
\end{enumerate}

One can adapt the foregoing to this case, more or less,
by treating $v_*$ as if it has degree $d$. 
But let us first isolate the cases not already
explicitly covered by our constructions.

First, suppose 
\begin{enumerate}
\item [$B_1$]
There is an external vertex of
degree $d$ which is not adjacent to $v_*$.
\end{enumerate}
Let $C$ be a $v_*$-component containing 
a vertex of degree $d$ not adjacent to $v_*$, 
and let $\hat C$ be the induced subgraph of $T$
on $C$ with the vertex $v_*$ adjoined.

If there are two nonadjacent vertices of degree $d$
in $\hat C$ then one of the foregoing cases applies.
Otherwise, considering the vertices of degree $d$
in $\hat C$ together with $v_*$ we have
one of the following possibilities:
\begin{enumerate}
\item[$B_1$.1] There is a unique vertex $v_1$
of degree $d$ in $\hat C$, 
at distance $\ell\ge 2$ from $v_*$
\item[$B_1$.2] There are two adjacent vertices $v_1,v_1'$
of degree $d$ in $\hat C$.
\end{enumerate}
In either case, if $v_*$ has degree $d$, with $d\ge 4$,
we will fall 
into one of our cases
with $v_0=v_*$. We need to adapt our
constructions when $v_*$ does not have degree $d$.

The effect of this is that in each case 
we need to reexamine the
proof that the resulting graph is $T$-free.

We will discuss these cases further below.

We have also the following possibility to consider.
\begin{enumerate}
\item [$B_2$]
$T'$ is a near-path with center $v_*$.
Every vertex of degree $d$ other than $v_*$
is adjacent to $v_*$.
\end{enumerate}
Here $v_*$ may or may not have degree $d$ itself.

This case escapes from those treated earlier and
must be handled separately. 

\subsection{The case $d=3$}

It will be convenient to clear away the case in which
the maximal degree $d$ is $3$. In particular in
this case the center $v_*$ has degree $d$.

If some external branch vertex has distance at least $2$ from the
closest branch vertex in $T$ then Case I or II applies.
So we may suppose that every external branch vertex
is adjacent to a branch vertex. 

If some $v_*$-component $C$ of $T$ contains at least two branch vertices
then taking $v_1'$ the external branch vertex of $C$, 
$v_1$ the neighboring branch vertex, and $v_0$ the next branch vertex
in $T$, that is either the next one in $C$ or $v_*$ itself,
we fall into Case IIIB, IIIC, or IIID.

So we may suppose that 
$$\hbox{All branch vertices other than $v_*$ are adjacent to
  $v_*$.}\leqno(*)$$ 

\begin{texteqn}[0.70]{Case IVA}
$T'$ is a near-path, $d=3$, and there are exactly two branch vertices 
in $T$, namely the center $v_*$, and an adjacent vertex
$v_1$.
\end{texteqn}

This case 
is highly reminiscent of Case IIIA, 
and we may adapt that construction
as follows. 

\begin{construction}
Take a two-way infinite path $P$ and divide it into intervals
$P_i,Q_i$ 
which are alternately of lengths $p_i=4$ and 
$q_i=1$ or $2$.
Adjoin a common neighbor $b_i$ to each interval $P_i$.
Call the result $\Gamma$.
\end{construction}

\begin{lemma}
$\Gamma$ is $T$-free
\end{lemma}
\begin{proof}
Consider any subgraph of $\Gamma$ containing two
adjacent branch vertices. These must be of the
form $u,b_i$ with $u\in P_i$. Suppose that this
graph is part of a subgraph of $\Gamma$ isomorphic with $T$.
Then either $u$ or $b_i$ corresponds to $v_*$
and hence is an endpoint of three disjoint paths 
of length $2$. The vertex $b_i$ as well as
any endpoint of $P_i$ has this property, but the three
paths involved must completely cover $P_i$ and $b_i$. 
So if the other vertex is to represent a branch point
of the subgraph, it cannot be $b_i$. Thus $b_i$ must
play the role of $v_*$ and then the three paths must embed
in such a way as to cover the neighbors of both endpoints of $P_i$ as
well and again $u$ cannot be a branch vertex of the image.
\end{proof}

\begin{lemma}
If $\Gamma\includedin G$ and $G$ is $T$-free,
then the graph induced on $V(\Gamma)$ by $G$ is 
a connected component of $G$, and consists of $\Gamma$ 
with any additional edges $(u,v)$ involving vertices $u,v\in P$ at
distance $2$, of the following types:
\begin{enumerate}
\item For some $q_i=2$, if $u\in Q_i$ is
adjacent to $u'\in P_{i\pm 1}$ 
then possibly $(u,v)$ is an edge, with $v$ the neighbor
of $u'$ in $P_{i\pm 1}$.
\item Two vertices of some $P_i$ at distance 2 (one
an endpoint, one an interior point) may be adjacent.
\end{enumerate}
\end{lemma}
\begin{proof}
First we exclude edges between $\Gamma$ and $G\setminus \Gamma$. 
If $(u,v)$ is an edge with $u\in \Gamma$ and
$v\notin \Gamma$ we look for an embedding
of $T$ into $\Gamma$ in which $u$ plays the role of $v_1$,
and we need to select an appropriate vertex to play the role of $v_*$.

If $u$ is adjacent to an endpoint $u'$ of $P_i$ then
$u'$ can play the role of 
$v_*$, though there are two cases to be distinguished here:
 $u=b_i$ or $u\in P$. We leave further inspection of this case  to the
reader.

Now suppose $u$ is not adjacent to an endpoint of $P_i$.
Then $u$ is an endpoint of $P_i$ itself, and $b_i$
can play the role of $v_*$.

So the graph $G_0$ induced on the vertices of $\Gamma$
by $G$ is a connected component of $G$ and the question remains as to
its precise structure.

Suppose first that $(b_i,v)$ is an edge in $G$
but not in $\Gamma$. Then easily $v\in P$ and since $v\notin P_i$ we
can
take $v$ to play the role of $v_*$ and $b_i$ to play the
role of $v_1$, taking as leaf adjacent to $b_i$ 
an interior vertex of $P_i$ lying on the side away from $v$.

So the only edges that come into consideration are edges
$(u,v)$ joining vertices of $P$. 
If $d(u,v)\ge 3$ and $u\in P_i$ 
then we let $u$ play the role of $v_1$ with 
$b_i$ representing a leaf adjacent to $u$, and $v$ may play the role
of $v_*$.

If $d(u,v)\ge 3$ and $u\in Q_i$ with $u$ adjacent to $u'\in P_{i\pm 1}$
then we let $u$ play the role of $v_*$ with
$u',b_i$ representing a path attached to $u$ and with
$v$ in the role of $v_1$.

So we have
$$u,v\in P; \quad d(u,v)=2$$

Now if $u\in Q_i$ and $q_i=1$ then easily as $u$ is adjacent to endpoints
in $P_{i-1}$ and $P_{i+1}$, there can be no new neighbors
of $u$ in $P$.

If $u\in Q_i$ with $q_i=2$ 
is adjacent to $u'$ in $P_{i+1}$
then there can be no edge $(u,v)$ with $v$ the
endpoint of $P_i$ closest to $u$, as then $v$ could play the role of
$v_1$ and
$b_i$ could play the role of $v_*$.

So for $u\in Q_i$ we have only the case mentioned in the statement of
the lemma.
\end{proof}

\begin{proposition}
In case IVA, there is no weakly universal $T$-free graph.
\end{proposition}
\begin{proof}
One needs to decode some information from an embedding
of $\Gamma$ into a $T$-free graph $G$.
This can be simplified by taking
$q_i=1$ over large intervals of fixed size,
with occasional values of $q_i=2$ inserted optionally
at regular intervals.

Let $G_0$ be the graph induced on $V(\Gamma)$ by $G$.
Viewing $G_0$ as a collection of $2$-connected blocks
which are connected in a tree structure, we see that the
$2$-connected blocks have approximately the same
vertices as they do in $\Gamma$, with some possible
variation involving $Q_i$ when $q_i=2$. The vertices
with four neighbors in their $2$-block are the $b_i$
and possibly some interior vertices of $P_i$. 

Most of the nontrivial 2-connected blocks
have order $5$,
and their points of attachment are the endpoints of the intervals $P_i$.
The exceptions may occur when $q_i=2$ and these occurrences will be
signalled either by the presence of a 2-connected block of order
$6$, or by two successive $2$-connected blocks with a gap of size
$2$. From this rudimentary analysis one cannot recover
the exact placement of the exceptional values of $q_i$,
but one can localize it with an error of $\pm 1$, which
is good enough.
\end{proof}

The next case to consider would have $T'$ a near-path,
exactly three branch vertices consisting of the center
$v_*$ and two neighbors of $v_*$, but Case IIIB covers
this one.

So in fact there is just one more case with $d=3$.

\begin{texteqn}[0.70]{Case IVB}
$T'$ is a near-path, $d=3$, 
and there are exactly four branch vertices 
in $T$, namely the center $v_*$, and three adjacent vertices
$v_1, v_2, v_3$, in distinct $v_*$-components.
We may suppose that $v_3$ is adjacent to two leaves
(and the same may possibly 
apply to one or both of $v_1,v_2$).
\end{texteqn}

\begin{construction}
Begin with a two-way infinite path $P$ divided into 
intervals $P_i,Q_i$ of lengths $p_i=6$ and $q_i=1$ or $0$. 
Attach a vertex $b_i$ adjacent to the vertices of $P_i$
and if $Q_i$ contains a vertex, give it infinite degree.
Call the result $\Gamma$.
\end{construction}

\begin{lemma}
$\Gamma$ is $T$-free
\end{lemma}
\begin{proof}
Consider the subgraph $\Gamma_0$ of $\Gamma$
with the same vertices, and with edges between
any pair of vertices $u,v$ of $\Gamma$
which lie in a subgraph of $\Gamma$ 
for which $u,v$ are adjacent branch vertices
of degree $3$ with no common neighbor.
Then all the edges containing $b_i$ in $\Gamma$
are retained, but the only edges along $P$ which
are retained are the ones involving a vertex of 
some $Q_i$, and the vertices of $Q_i$ have degree $2$
in $\Gamma_0$.

Under an embedding of $T$ into $\Gamma$, $v_*$ must correspond 
to a vertex of degree at least $3$ in $\Gamma_0$, 
thus a vertex $b_i$ for some $i$. But there is no such embedding 
as one of the three branch vertices neighboring $v_*$
must correspond to an interior point of $P_i$, 
so that together with its neighbors this $v_*$-component
requires at least $3$ vertices of $P_i$, and each of the
others requires at least $2$ vertices of $P_i$.
\end{proof}

\begin{lemma}
If $\Gamma$ is contained in the $T$-free graph $G$,
and if $G_0,\Gamma_0$ are the subgraphs of $G$
induced on the branch vertices of $G$ and of $\Gamma$ respectively, 
then $\Gamma_0$ is a connected component of $G_0$.
\end{lemma}
\begin{proof}
The vertices of $\Gamma_0$ are the vertices of $P$
together with the vertices $b_i$.

We claim first
$$\hbox{There is no edge $(u,v)$ in $G$ with $u\in P_i$,
$v\notin \Gamma_0$}\leqno(1)$$
We may suppose that $u$ is the first, second, or third
vertex of $P_i$. If $u$ is the first or third
vertex then we embed $T$ into $\Gamma$ with
$b_i$ playing the role of $v_*$ and with the extra
vertex $v$ playing the role of a leaf attached
to one of its neighbors.
If $u$ is the second vertex of $P_i$ we let
the first vertex of $P_i$ play the role of $v_*$.
So $(1)$ holds.

Now we claim 
\begin{texteqn}{(2)}
There is no edge $(u,v)$ in $G$ with
$u\in Q_i$, $v\notin \Gamma_0$, and 
$v$ a branch vertex of $G$
\end{texteqn}
Let $u,u',u''$ be three neighbors of $v$ in $G$.
We attempt to embed $T$ into $G$ with $u$
playing the role of $v_*$. This can be blocked if
$u'$ or $u''$ lies on $P$ or coincides with $b_i$ or 
$b_{i+1}$. 
Suppose therefore that $u'$ is of one of these two forms.

If $u'=b_{i+1}$ then we may let $b_{i+1}$ play the role of 
$v_*$
with $u'$ as one of its neighboring branch vertices.
This could only be blocked by having
$u''\in P$ in which case $u''\in Q_j$ for some $j$. 
In this case $u''$ could play the role of $v_*$
instead.
Similarly the case $u'=b_i$ may be excluded.

So we may suppose that the vertices $u',u''$ which lie
in $\Gamma$ lie in $P$, and in each such case in some $Q_j$. 
Furthermore we may suppose that among those vertices
of $u,u',u''$ lying on $P$, $u$ is the first in order.
Then we use $u$ to represent $v_*$, getting a contradiction. 
So $(2)$ holds.

So our lemma is proved as far as edges involving vertices of $P$ are concerned.
Suppose finally that $(u,v)$ is an edge with $u=b_i$
for some $i$ and $v$ a branch vertex of $G$ not in 
$\Gamma_0$. Let $u',u''$ be additional neighbors of $v$.
By the cases already treated, these vertices do
not lie on $P$. It is then easy to embed $T$ into $G$
with $b_i$ playing the role of $v_*$, arriving at a contradiction.
\end{proof}

\begin{lemma}
If $\Gamma$ is contained in the $T$-free graph $G$,
and $\Gamma_0$ is the graph induced on the branch vertices
of $\Gamma$ by $G$, then any edge $(u,v)$ of $\Gamma_0$
which is not an edge of $\Gamma$ involves two vertices 
of $P$, at distance at most $3$, and of one of the following two
forms. 
\begin{enumerate}
\item $u$ or $v$ is in $Q_i$ for some $i$;
\item $u$ and $v$ are in successive intervals
$P_{i-1},P_i$, with $q_i=0$, and adjacent 
to endpoints of these intervals; $d(u,v)=3$.
\end{enumerate}
\end{lemma}
\begin{proof}
First, one may eliminate the possibility $u=b_i$,
as a new edge of this type leads to an embedding
of $T$ into $G$ with $b_i$ playing the role of $v_*$. 
So we may suppose $u,v\in P$.

Now suppose $d(u,v)\ge 4$ where all distances will
be measured in $P$.
Then we let $u$ play the role of $v_*$ and
we let the neighbor $u'$ of $u$ along $P$ in the direction of $v$ play
the role of a branch vertex adjacent to $v_*$, whose
further neighbors are leaves of $T$.
As $d(u,v)\ge 4$ this embedding can be completed to
an embedding of $T$.

If $u$ is an endpoint of $P_i$, say a left endpoint,
and $v$ lies farther to the left along $P$,
embed $T$ into $G$ with $b_i$ representing
$v_*$ and with the immediate neighbor of $u$
to its left representing a leaf adjacent to $u$.

If $v$ lies to the right of $P$, and at distance at most $3$,
then let $u$ play the role of $v_*$ with $v$ an adjacent branch
vertex. Here $b_i$ will also play the role of an adjacent branch
vertex.

In the remaining cases, we may choose notation
so that $u$ is adjacent to an endpoint $u'$ of $P_i$
for some $i$. Leaving aside the cases mentioned
in the statement of the lemma, we may suppose
$v\in P_i$ as well. We let $u'$ represent $v_*$
and $u$ represents a branch vertex adjacent to $u'$
with a neighbor on $P$ and $v$ as its adjacent leaves. 
\end{proof}

\begin{proposition}
In case IVB, there is no weakly universal $T$-free graph.
\end{proposition}
\begin{proof}
If $\Gamma$ is contained in a $T$-free graph $G$
then the graph $\Gamma_0$ induced on the branch vertices
of $\Gamma$ by $G$ can be recovered from one of its vertices.

We examine the vertices of degree $6$ in 
$\Gamma_0$. These include the vertices $b_i$, and for
these vertices the graph induced on its neighbors is connected. 

If $q_i=1$ and $u$ is the unique vertex of $Q_i$,
then $u$ can have a maximum of 6 neighbors in $\Gamma_0$,
which would then be all of its neighbors in $P$
up to distance $3$. In this case the graph
induced on the neighbors of $u$ in $\Gamma_0$ is disconnected. 

If $u$ lies in some interval $P_i$,
then in addition to its $3$ neighbors in $\Gamma$,
there can be at most one more in $\Gamma_0$.
So these do not come into consideration.

As we may distinguish the $b_i$ by the structure of
the graph induced on their neighbors,
we can also recognize the path $P$ and the 
intervals $P_i$, which remain paths in $\Gamma_0$.
While $\Gamma_0$ may have some additional edges
we can then detect the vertices in $Q_j$
for $q_j=1$, as well as their locations relative
to the $P_i$.
\end{proof}

\subsection{The case $d\ge 4$}

$T'$ is a near-path with center $v_*$. 
If the degree of $v_*$ is $d$, or if one of
the $v_*$-components of $T$ contains two vertices
of degree $d$, then one of the cases I, II applies.

So we suppose 
\begin{enumerate}
\item $\deg(v_*)<d$.
\item Each $v_*$-component of $T$ contains
at most one vertex of degree $d$.
\end{enumerate}

As we have disposed of the case in which there is a unique vertex of
degree $d$ in $T$, there are either two or three vertices
of degree $d$.

\begin{texteqn}[0.70]{Case IVC}
$T'$ is a near-path, $d\ge 4$, and there are exactly
two vertices $v_0,v_1$ 
of degree $d$ in $T$, lying in distinct
$v_*$-components.
\end{texteqn}

If $d(v_0,v_1)=2$ then Case IC applies, so 
we suppose $\ell=d(v_0,v_1)\ge 3$.

\begin{construction}
Take a $(d-1)$-regular tree $P$
and find disjoint intervals $P_i$ in $P$ 
of lengths $p_i\ge 3\ell-3$ 
such that every vertex of $P$ not in one of the $P_i$
lies at distance
less than $\ell$ from at least two of the intervals $P_i$, 
but no vertices of distinct $P_i$, $P_j$ lie at distance less than
$\ell$ of each other.

This is done inductively. 
At each stage finitely many intervals $P_i$
have been selected, no two at distance less than $\ell$,
so that the subgraph induced on the vertices within distance $\ell-1$
of some $P_i$ is connected. 
Then a vertex at minimal distance $\ell$ from some
$P_i$ is selected and put into a new interval $P_i$.
With some housekeeping one may ensure that the
whole tree $P$ is exhausted by this process.
For any vertex $v$ of $P$ which is not in the $P_i$, there
is a first stage at which $v$
falls within distance $\ell-1$ of one of the $P_i$,
and at that stage $v$ lies outside the convex
hull of the $P_i$ selected up to that point.
Consider the tree $T_v$ rooted at $v$ obtained
by deleting the $v$-component of $P$ containing the
convex hull of the $P_i$. 
If the distance from $v$ to the nearest $P_i$ (so far chosen)
is $k$, then the vertices of $T_v$ lying at distance $\ell$
from the closest $P_i$ are those at distance $\ell-k$ from $v$.
As the construction proceeds, and new intervals $P_i$ are selected,
one of two things will occur.
Possibly the distance from $v$ to the nearest $P_i$ will
be diminished at some point, in which case $v$ is close
to at least two such intervals. In the contrary case,
the distances of the vertices in $T_v$ to the nearest
$P_i$ will also be unaltered until
one of them, lying at distance $\ell$ from the
nearest $P_i$ and at distance $\ell-k$ from $v$, 
is selected as an endpoint of a new
interval $P_i$, at which point $v$ will lie within $\ell-1$ of at
least two such intervals.

Now adjoin a vertex $b_i$ adjacent to the vertices
of $P_i$ for each $i$, and attach to $b_i$
infinitely many $(d-1)$-regular trees---but with its root
of degree $(d-2)$, so that after attachment to $b_i$ 
the degree of the root is also $(d-1)$.

Call the result $\Gamma$.
\end{construction}

\begin{lemma}
The graph $\Gamma$ is $T$-free.
\end{lemma}
\begin{proof}
The only vertices of degree $d$ in $\Gamma$
are the $b_i$ and the vertices in the intervals $P_i$.
As these are widely spaced, the only way to embed
$T$ into $\Gamma$ is to send the two vertices of degree
$d$ into $P_i\union \{b_i\}$ for some $i$.

The vertices of $P_i$ have degree $d$ and a common neighbor, so at
most one of these vertices can serve as the image of a vertex of
degree $d$ in $T$. Therefore the two images must be of the form
$b_i,u$ with $u\in P_i$. But then the edge $(u,b_i)$ is 
not in the image, and $u$ cannot have degree $d$ in the image.
\end{proof}

\begin{lemma}
If $\Gamma$ is contained in a $T$-free graph $G$,
then 
\begin{enumerate}
\item The vertices of $P$ have finite degree in $G$,
and their neighbors all lie in $\Gamma$.
\item For $u\in P$, the neighbors of $u$ of finite
degree in $G$ are its neighbors in $P$.
\end{enumerate}
\end{lemma}
\begin{proof}
First, if one adjoins an edge $(u,v)$ linking
a vertex $u\in P$ to
a vertex $v\notin \Gamma$, then $T$ embeds into the extended graph by
finding a path of length $\ell$ along $P$ to some $b_i$.
Note that in view of the structure of $T$ the vertices of degree $d$
have at least one adjacent leaf, so the new vertex $v$
can serve to represent one such leaf, $u$ can represent a vertex of
degree $d$, and $b_i$ can represent the other vertex of degree $d$.
As there are only two vertices of degree $d$ in $T$,
the embedding may be completed.

The same applies if the vertex $v$ lies in one
of the trees attached to a $b_i$. So we may
suppose that $v$ lies on $P$ or among the vertices
$b_i$, and not too far from $u$. So the first point follows. 

If we now require $v$ to have finite degree then
we are no longer concerned with the $b_i$,
and we may suppose $v\in P$.

Now if $v$ is not adjacent to $u$ in $P$ there are
various possibilities. Let us fix a leaf $\tilde v$
adjacent to $v_0$ and an embedding of $T\setminus \tilde v$ into 
$\Gamma$ taking one
vertex of degree $d$ to $u$ and the other 
to some $b_i$. 
We can extend this to an embedding of $T$ into $G$
unless $v$ lies either on the path
from $u$ to $b_i$ along $P$, or in the remaining part
of the neighborhood of $u$ in $\Gamma$ used to embed
the other $v_0$-components of $T$.

In the second case, we can examine the path from $u$ to $v$ along
$P$ and use the neighbor $u'$ of $u$ along this path to represent 
$\tilde v$, and the tree originating with the edge from $u$ through $v$ and
continuing along $P$ to replace the $u$-component of $P$ 
containing $u'$.

So suppose that $v$ lies along the path $L$ from $u$
to a neighbor of $b_i$ in $P_i$, of length $\ell-1$.
One or both of the vertices $u,v$ may lie in the interval $P_i$. 
Removing such vertices, $P_i$ is divided into at most
three intervals, of total length at least $3\ell-5$,
and hence one of these intervals contains
at least $\ell-1$ vertices.
Such an interval may or may not be separated from $u$
by $v$. If it is not separated, we can make use
of it and possibly other vertices of $P$ to find a path
of length $\ell$ from $u$ to $b_i$ avoiding $v$, and
complete the construction. If it is separated, we
can make use of the edge $(u,v)$ to find a suitable
replacement path, and use the neighbor of $u$ on $P$ in the direction
of $v$ as a representative for the leaf $\tilde v$.
\end{proof}

\begin{proposition}
In case IVC, there is no weakly universal $T$-free graph.
\end{proposition}
\begin{proof}
We have considerable latitude in the choice of 
the size $p_i$ of $P_i$.
It suffices to decode the set of $p_i$ involved in the
construction after $\Gamma$ is embedded into a $T$-free
graph $G$ (whereas 
the ``sequence'' is not that well-defined
at this point).

By the preceding lemma, we can recover the graph
structure on $P$ from one of its vertices, in $G$.
We would like to recover the intervals $P_i$
by considering vertices of $P$ with a common neighbor
in $G$ lying at the root of an infinite system 
of $(d-1)$-regular trees.
By the preceding lemma this common neighbor
would have to lie in $\Gamma$ and 
it only be some $b_i$.
So we have to deal with the possibility that a vertex
$u$ outside the interval $P_i$ might be
connected to $b_i$. But then by our construction, this vertex 
would lie
within $(\ell-1)$ of some second interval $P_j$, 
leading to an embedding of $T$ into $G$. 
\end{proof}

\medskip  

Finally we come to the case of three vertices
of maximal degree, with the center $v_*$ 
of lower degree. 
At least one of these three vertices must be adjacent to
$v_*$. We can unify the treatment of these cases, but
we prefer to first treat the case in which all vertices of degree $d$
are adjacent to $v_*$, and then discuss the modification of our
construction suitable for 
other cases.

\begin{texteqn}[0.70]{Case IVD}
$T'$ is a near-path, $d\ge 4$, and 
there are three vertices $v_0,v_1,v_2$
of degree $d$ in $T$,  all of which are adjacent to $v_*$.
\end{texteqn}

\begin{construction}
Let $C_1,C_2$ be the $v_*$-components of $T$ containing
$v_1,v_2$ respectively.
Let $H$ be the graph obtained by amalgamating the induced graph on  
$\{v_*\}\union C_1\union C_2$ 
with itself over the vertices
$v_*,v_1,v_2$ and then giving $v_1,v_2,v_*$ 
infinitely many new neighbors.

Take a two-way infinite path $P$ partitioned into intervals $P_i$
of lengths $p_i=1$ or $2$. Adjoin a vertex $b_i$
adjacent to each vertex of $P_i$, and attach
$H$ to $b_i$ with $b_i$ playing the role of $v_*$.

Raise the degrees of the vertices on $P$ to
$d-1$ by adding additional edges between pair on $P$,
keeping the girth of the induced graph on $P$ extremely large
(it will resemble a $(d-2)$-regular tree locally).

Then take any remaining vertices of degree less than $d-1$ and
attach trees to them so as to raise all such vertex 
degrees up
to $d-1$.

Call the result $\Gamma$.
\end{construction}

\begin{lemma} The graph $\Gamma$ is $T$-free.
\end{lemma}
\begin{proof}
There are no vertices with three neighbors of degree $d$.
\end{proof}

\begin{lemma}
If $\Gamma$ is contained in a $T$-free graph $G$ 
then the vertices of $P$ have finite order in $G$,
and if $u\in P$ then the neighbors $v$ of $u$ of finite order in $G$
are its neighbors in $P$.
\end{lemma}
\begin{proof}
A vertex $u\in P$ can have no new neighbor $v$ in 
$G\setminus \Gamma$ as this immediately produces an embedding of $T$
into $\Gamma$. 

If $u\in P_i$ then the only candidates for a vertex $v\in \Gamma$
which could serve as a new neighbor without producing an
embedding of $T$ into $\Gamma$ are the vertices
adjacent to $b_i$ which correspond to $v_1$ or $v_2$, and these have
infinite degree in $\Gamma$.
So the vertices of $P$ have finite degree,
and if we restrict our attention to neighbors
of finite degree then there are none available
other than the $(d-1)$ neighbors we have already selected.
\end{proof}

\begin{proposition}
In case IVD, there is no weakly universal $T$-free graph.
\end{proposition}
\begin{proof}
If $\Gamma$ is contained in the $T$-free graph $G$
then we can recover the graph $P$ from $G$
and one vertex of $P$.
Now if $d\ge 5$ there is enough variability in the
structure of $P$ itself to yield the desired conclusion,
so we may suppose that $d=4$ and $P$ is a path.

If a pair $(u,u')$ 
of vertices of $P$ has a common neighbor,
then that neighbor is not on $P$ and could only be
some $b_i$ or some vertex adjacent to $b_i$,
and then only if $u,u'\in P_i$. 
Thus the sets $P_i$ can be recovered, and thus 
the sequence $p_i$ can be recovered up to a shift
and reversal.
\end{proof}

\begin{texteqn}[0.68]{Case IVD$'$}
$T'$ is a near-path, $d\ge 4$, and 
there are three vertices $v_0,v_1,v_2$
of degree $d$ in $T$, in distinct $v_*$-components of $T$.
\end{texteqn}

\begin{construction}
We proceed much as in the previous case but
with a different treatment for the $v_*$-components
$C_i$ containing $v_i$ nonadjacent to $v_*$ ($i=1$ or $2$, possibly): 
these components we allow to be freely amalgamated
over $v_*$ (without fixing the vertex $v_i$).
Otherwise, we proceed as in the previous construction.

Call the result $\Gamma$.
\end{construction}

\begin{lemma}
The graph $\Gamma$ is $T$-free.
\end{lemma}
\begin{proof}
If $v_*$ is adjacent to $n$ vertices of degree $d$
in $T$, where $1\le n\le 3$,
then no vertex of $\Gamma$ is adjacent to
more than $n-1$ vertices of degree $d$.
\end{proof}

This is actually the main point, since we have loosened the
construction of $H$ in a way that in other contexts
could easily lead to a violation of this first step.

The rest of the analysis is as before since the construction
of $H$ is if anything even freer than it was in the previous case.

\begin{proposition}
In case IVD$'$, there is no weakly universal $T$-free graph.
\end{proposition}

\medskip
With this, the proof of Theorem \ref{TC1:main}$'$ is complete, and
thus also the proof of Theorem \ref{TC1:main}. Together with Theorem 
\ref{Monarchy:main} this gives the full Tree Conjecture.

We believe that these methods can be applied to the complete
identification of all finite connected graphs $C$ for which there is a
countable universal (weakly or strongly) $C$-free graph, in part because we
expect the list of exceptional $C$ allowing such a universal graph to
be fairly limited. 
In particular we now think it quite likely that this problem is decidable
in the case of a single constraint, and very possibly more generally.

Based on the results of \cite{ct} it would appear
that the ``generic'' case corresponds roughly to the case in which
there is some block of order at least 6, and that the nongeneric case 
is therefore inconveniently complicated.

We remark that a proof of decidability for the case of a single
constraint may be achievable without
actually working through all the critical cases. Once the set of
unsolved cases is reduced to a well-quasiordered set relative 
to the ``pruning'' relation, one knows that the remaining minimal cases not
allowing a universal graph form a finite set, and the problem is therefore 
algorithmically decidable. 
This style of argument does not necessarily 
provide any further indication as to
what the relevant finite subset might be, any bound on its size, or a
fortiori any explicit algorithms.
We hope to return to this topic.

\section{List of cases}

\subsection{Paths and some near-paths}

\begin{texteqn}[0.70]{Case I}
$\ell \ge 2$
\end{texteqn}

\begin{texteqn}[0.70]{Case IA}
$T$ has a vertex $v_0$ of maximal degree $d\ge 4$ such that some
  $v_0$-component $C$ of $T$ contains a unique vertex $v_1$ of degree
  $d$, and  $\ell=d(v_0,v_1)\ge 3$. 
Either $T'$ is a path, or else $T'$ is a near-path
whose center does not lie in the $v_0$-component $C$. 
\end{texteqn}

\begin{texteqn}[0.70]{Case IB}
$T$ has a vertex $v_0$ of maximal degree $d=3$ such that some
  $v_0$-component $C$ of $T$ contains a unique vertex $v_1$ of degree
  $d$, and  $\ell=d(v_0,v_1)\ge 3$. 
Either $T'$ is a path, or else $T'$ is a near-path
whose center does not lie in the $v_0$-component $C$. 
\end{texteqn}

\begin{texteqn}[0.70]{Case IC}
$T$ has a vertex $v_0$ of maximal degree $d$ such that some
  $v_0$-component $C$ of $T$ contains a unique vertex $v_1$ of degree
  $d$, and  $\ell=d(v_0,v_1)= 2$. 
Either $T'$ is a path, or else $T'$ is a near-path
whose center is not $v_1$.
\end{texteqn}

\begin{texteqn}[0.70]{Case II}
$\ell =1$
\end{texteqn}

\begin{texteqn}[0.70]{Case IIA}
$T$ has a vertex $v_0$ of maximal degree $d\ge 5$ such that some
$v_0$-component $C$ of $T$ contains a unique vertex 
$v_1$ of degree $d$, and $v_0,v_1$ are adjacent.
Either $T'$ is a path, or $T'$ is a near-path whose center 
does not lie in the $v_0$-component $C$.
\end{texteqn}
\begin{texteqn}[0.70]{Case IIB}
$T$ has a vertex $v_0$ of maximal degree $d=4$ such that some
$v_0$-component $C$ of $T$ contains a unique vertex 
$v_1$ of degree $d$, and $v_0,v_1$ are adjacent.
Either $T'$ is a path, or $T'$ is a near-path whose center 
does not lie in the $v_0$-component $C$.
\end{texteqn}

\begin{texteqn}[0.70]{Case IIIA}
$T'$ is a path and $T$ contains exactly two branch vertices $v_0$ and  
$v_1$, which are adjacent, and of degree $3$.
\end{texteqn}

\begin{texteqn}[0.70]{Case IIIB}
The maximal vertex degree is $3$. $T$ contains a sequence of three
adjacent branch vertices $v_1,v_0,v_2$ with $v_1$ external and
adjacent to a leaf. Some $v_0$-component of $T$ is a path. 
\end{texteqn}

\begin{texteqn}[0.70]{Case IIIC}
The maximal vertex degree is $3$. 
$T$ contains a sequence of three successive
branch vertices $v_0,v_1,v_1'$ with $v_1'$ external and
adjacent to a leaf, $v_1$ adjacent to $v_1'$,
and $\ell=d(v_0,v_1)\ge 3$, 
where $v_0$ is the closest branch vertex to
$v_1$ other than $v_1'$. 
Either the pruned tree $T'$ is a path, 
with both external branch vertices of
$T$ are adjacent to branch vertices. or a near-path
with the center not in the $v_0$-component containing $v_1$. 
\end{texteqn}

\begin{texteqn}[0.70]{Case IIID}
The maximal vertex degree is $3$. 
$T$ contains a sequence of three successive
branch vertices $v_0,v_1,v_1'$ with $v_1'$ external and
adjacent to a leaf, $v_1$ adjacent to $v_1'$,
and $\ell=d(v_0,v_1)=2$, 
where $v_0$ is the closest branch vertex to
$v_1$ other than $v_1'$. 
Either the pruned tree $T'$ is a path, 
and both external branch vertices of
$T$ are adjacent to branch vertices, or it is a near-path
with the center not in the $v_0$-component containing $v_1$.  
\end{texteqn}

\subsection{The remaining near-paths}

\begin{texteqn}[0.70]{Case IVA}
$T'$ is a near-path, $d=3$, and there are exactly two branch vertices 
in $T$, namely the center $v_*$, and an adjacent vertex
$v_1$.
\end{texteqn}

\begin{texteqn}[0.70]{Case IVB}
$T'$ is a near-path, $d=3$, 
and there are exactly four branch vertices 
in $T$, namely the center $v_*$, and three adjacent vertices
$v_1, v_2, v_3$, in distinct $v_*$-components.
We may suppose that $v_3$ is adjacent to two leaves
(and the same may possibly 
apply to one or both of $v_1,v_2$).
\end{texteqn}

\begin{texteqn}[0.70]{Case IVC}
$T'$ is a near-path, $d\ge 4$, and there are exactly
two vertices $v_0,v_1$ 
of degree $d$ in $T$, lying in distinct
$v_*$-components.
\end{texteqn}

\begin{texteqn}[0.70]{Case IVD}
$T'$ is a near-path, $d\ge 4$, and 
there are three vertices $v_0,v_1,v_2$
of degree $d$ in $T$, and all are adjacent to $v_*$.
\end{texteqn}

\begin{texteqn}[0.68]{Case IVD$'$}
$T'$ is a near-path, $d\ge 4$, and 
there are three vertices $v_0,v_1,v_2$
of degree $d$ in $T$, in distinct $v_*$-components of $T$.
\end{texteqn}




\begin{thebibliography}{10}


\bibitem{ck}
G.~Cherlin and P.~Komj\'{a}th
\newblock  There is no universal countable pentagon-free graph. 
\newblock {\sl J.~Graph Theory}, {\bf 18} (1994), 337--341.

\bibitem{css}
G.~Cherlin, S.~Shelah, and N.~Shi
\newblock  Universal graphs with forbidden subgraphs and algebraic closure
\newblock Advances in Applied Mathematics {\bf 22} (19??), ??-??

\bibitem{cs1}
G.~Cherlin and N.~Shi. 
\newblock Graphs omitting a finite set of 
cycles.
\newblock {\sl J.~Graph Theory}, {\bf 21} (1996), 351--355.

\bibitem{cs2}
G.~Cherlin and N.~Shi. 
\newblock Graphs omitting sums of complete graphs, 
\newblock {\sl J.~Graph Theory}, {\bf 24} (1997), 237--247.

\bibitem{cs3}
G.~Cherlin and N.~Shi. 
\newblock Forbidden subgraphs and forbidden substructures.
\newblock {\sl J.~Symbolic Logic}, {\bf 66} (2001), 1342--1352.

\bibitem{cst}
G.~Cherlin, N.~Shi, and L.~Tallgren, 
\newblock Graphs omitting a bushy tree.
\newblock {\sl J.~Graph Theory}, {\bf 26} (1997), 203--210.

\bibitem{ct}
G.~Cherlin and L.~Tallgren, 
\newblock Graphs omitting a near-path or 2-bouquet.
\newblock Submitted.

\bibitem{fk1}
Z.~F\"{u}redi and P.~Komj\'{a}th. 
\newblock On the existence of countable universal graphs. 
\newblock {\sl J.~Graph Theory}, {\bf 25} (1997), 53--58.

\bibitem{fk2}
Z.~F\"{u}redi and P.~Komj\'{a}th. 
\newblock Nonexistence of universal graphs without some trees, 
\newblock {\sl Combinatorica}, {\bf 17} (1997), 163--171.

\bibitem{gk}
M.~Goldstern and M.~Kojman. 
\newblock Universal arrow free graphs, 
\newblock {\sl Acta Math.~Hungary} {\bf 73} (1996), 319--326.

\bibitem{kom}
P.~Komj\'{a}th.
\newblock Some remarks on universal graphs.
\newblock {\sl Discrete Math.} {\bf 199} (199), 259--265.

\bibitem{kmp}
P.~Komj\'{a}th, A.~Mekler and J.~Pach. 
\newblock Some universal graphs, 
\newblock {\sl Israel J.~Math}. {\bf 64} (1988), 158--168.

\bibitem{kp1}
P.~Komj\'{a}th and J.~Pach. 
\newblock Universal graphs without large bipartite subgraphs. 
\newblock {\sl Mathematika}, {\bf 31} (1984), 282--290.

\bibitem{kp2}
P.~Komj\'{a}th and J.~Pach. 
\newblock 
Universal elements and the complexity of certain classes of infinite graphs.
\newblock {\sl Discrete Math}. {\bf 95} (1991), 255--270.

\bibitem{pac}
J.~Pach. 
\newblock A problem of Ulam on planar graphs. 
\newblock {\sl Eur.~J.~Comb.~}{\bf 2} (1981), 357--361.

\bibitem{rad}
R.~Rado.
\newblock Universal graphs and universal functions.
\newblock {\sl Acta Arith.}~{\bf 9} (1964), 331--340.

\end{thebibliography}
\end{document}